\documentclass[oneside, reqno]{amsart}

\usepackage{xypic}
\input xy
\xyoption{all}
\usepackage{epsfig}
\usepackage{amsthm}
\usepackage{amssymb}
\usepackage{dsfont}
\usepackage{amsmath}
\usepackage{amscd}
\usepackage{color}
\usepackage[T1]{fontenc}
\usepackage{stackengine}
\usepackage{upgreek}
\usepackage[font=scriptsize]{caption}
\usepackage{wrapfig}
\usepackage[utf8x]{inputenc}
\stackMath

\newcommand{\bg}{\begin{equation}}
\newcommand{\ed}{\end{equation}}
\newcommand{\bga}{\begin{eqnarray}}
\newcommand{\eda}{\end{eqnarray}}
\newcommand{\pf}{\textbf{Proof:\ }}

\def\cbdu{\par{\raggedleft$\Box$\par}}

\newtheorem {Theorem}  {Theorem}

\numberwithin{Theorem}{section}

\newtheorem {Lemma}[Theorem]  {Lemma}
\newtheorem {Proposition}[Theorem]{Proposition}
\theoremstyle{definition}
\newtheorem{Definition}[Theorem]{Definition}
\theoremstyle{remark}

\newtheorem {Corollary}[Theorem]{\bf Corollary}

\expandafter\chardef\csname pre amssym.def
at\endcsname=\the\catcode`\@ \catcode`\@=11
\def\undefine#1{\let#1\undefined}
\def\newsymbol#1#2#3#4#5{\let\next@\relax
 \ifnum#2=\@ne\let\next@\msafam@\else
 \ifnum#2=\tw@\let\next@\msbfam@\fi\fi
 \mathchardef#1="#3\next@#4#5}
\def\mathhexbox@#1#2#3{\relax
 \ifmmode\mathpalette{}{\m@th\mathchar"#1#2#3}%
 \else\leavevmode\hbox{$\m@th\mathchar"#1#2#3$}\fi}
\def\hexnumber@#1{\ifcase#1 0\or 1\or 2\or 3\or 4\or 5\or 6\or 7\or 8\or
 9\or A\or B\or C\or D\or E\or F\fi}

\font\teneufm=eufm10 \font\seveneufm=eufm7 \font\fiveeufm=eufm5
\newfam\eufmfam
\textfont\eufmfam=\teneufm \scriptfont\eufmfam=\seveneufm
\scriptscriptfont\eufmfam=\fiveeufm

\catcode`\@=\csname pre amssym.def at\endcsname

\newcounter{remark}
\newenvironment{remark}
{\medskip \stepcounter{remark} \noindent \textit{Remark
\arabic{section}.\arabic{remark}.}}{\rm \cbdu}

\newcommand{\e}{\epsilon}

\newcommand{\om}{\omega}

\newcommand{\Om}{\Omega}
\renewcommand{\a}{\alpha}

\renewcommand{\b}{\beta}

\newcommand{\gm}{\gamma}

\def  \12  {{\frac{1}{2}}}

\def\build#1_#2^#3{\mathrel{\mathop{\kern 0pt#1}\limits_{#2}^{#3}}}

\numberwithin{equation}{section}

\begin{document}

\title[Stochastic Dyadic MHD model]{Uniqueness for a stochastic ideal dyadic MHD model}

\author [Mimi Dai]{Mimi Dai}

\address{Department of Mathematics, Statistics and Computer Science, University of Illinois at Chicago, Chicago, IL 60607, USA}
\email{mdai@uic.edu} 

\author [Qirui Peng]{Qirui Peng}

\address{Department of Mathematics, Statistics and Computer Science, University of Illinois at Chicago, Chicago, IL 60607, USA}
\email{qpeng9@uic.edu} 

\author [Cheng Ouyang]{Cheng Ouyang}

\address{Department of Mathematics, Statistics and Computer Science, University of Illinois at Chicago, Chicago, IL 60607, USA}
\email{couyang@uic.edu}


\begin{abstract}
We study a stochastic dyadic model with both forward and backward energy cascade mechanisms for the inviscid and non-resistive magnetohydrodynamics. For a particular class of stochastic forcing, we show weak uniqueness for the stochastic system. However the solution dissipates the energy which is formally an invariant quantity for the system.

\bigskip


\hspace{0.02cm}CLASSIFICATION CODE: 35Q35, 76B03, 76W05.
\end{abstract}

\maketitle

\section{Introduction}

\subsection{Dyadic models for MHD}
The incompressible inviscid and non-resistive (ideal) magnetohydrodynamics (MHD) governed by the equations
\begin{subequations}
\begin{align}
u_t+(u\cdot\nabla) u-(B\cdot\nabla) B+\nabla p=&\ 0, \label{mhda}\\
B_t+(u\cdot\nabla) B-(B\cdot\nabla) u =&\ 0, \label{mhdb}\\
\nabla \cdot u= &\ 0. \label{mhdc}
\end{align}
\end{subequations}
describes electrically conducting fluids in astrophysics in the context that the underlying length scales are large, and hence the kinetic viscosity and magnetic diffusivity are insignificant.  The unknowns are the fluid velocity $u$, magnetic field $B$, and pressure function $p$. As for the pure inviscid fluid featured by the Euler equation, i.e. (\ref{mhda}) with $B\equiv0$, understanding the dynamics of the ideal MHD thoroughly remains a great challenge. 
Inspired by the study of dyadic models for hydrodynamics (see \cite{CDF} and references therein), 
dyadic models were introduced in \cite{Dai-20} for MHD systems. 
In this paper, we consider a dyadic model for system (\ref{mhda})-(\ref{mhdc}) with both forward and backward energy cascade mechanisms,
\begin{equation}\label{sys-1}
\begin{split}
\frac{d}{dt}a_j=&-\left(\lambda_j^{\theta}a_ja_{j+1}-\lambda_{j-1}^{\theta}a_{j-1}^2\right)+\left(\lambda_{j}^{\theta}b_jb_{j+1}-\lambda_{j-1}^{\theta}b_{j-1}^2\right),\\
\frac{d}{dt}b_j= & -\left(\lambda_j^{\theta}a_jb_{j+1}-\lambda_{j}^{\theta}b_ja_{j+1}\right),
\end{split}
\end{equation}
for $j\geq0$ with $a_{-1}=b_{-1}=0$, and $\theta$ is a parameter associated with the intermittency effect of the dynamics. 
In the infinite ODE system (\ref{sys-1}), we identify the quantities $a_j=\|u_j\|_{L^2}$ and $b_j=\|B_j\|_{L^2}$, where $u_j$ and $B_j$ are the $j$-th shell Littlewood-Paley projections of $u$ and $B$ respectively.  
We observe that only interactions with the nearest neighbor shells are taken into account in the modeling. 
The model (\ref{sys-1}) preserves the most essential feature of the original dynamics 
of (\ref{mhda})-(\ref{mhdc}): the total energy 
\begin{equation}\notag
E(t)=\frac12\sum_{j\geq 0}\left(a_j^2(t)+b_j^2(t)\right)
\end{equation}
and the cross helicity defined by
\[H^c(t)=\sum_{j=0}^\infty a_jb_j\]
are conserved quantities formally.


If we take $b_j\equiv0$ for $j\geq0$, system (\ref{sys-1}) reduces to the dyadic Euler model, which has been well studied for decades. Although with similarly nonlinear structures with the dyadic Euler model, the MHD model (\ref{sys-1}) presents additional obstacles from the interactions of the fluid and magnetic field. In \cite{DF1}, the authors showed linear instability for (\ref{sys-1}) around certain steady state, although an analysis of nonlinear instability is out of reach due to the very reason of aforementioned intricate interactions. Nonetheless, for the dyadic MHD model with dissipation and external forcing, the authors \cite{DF2} constructed non-unique Leray-Hopf solutions in spite of the complicated nonlinear interactions.   

For the Euler dyadic model
\begin{equation}\label{euler}
\frac{d}{dt}a_j=-\left(\lambda_j^{\theta}a_ja_{j+1}-\lambda_{j-1}^{\theta}a_{j-1}^2\right),
\end{equation}
the authors of \cite{BFM} showed non-uniqueness of weak solutions in $l^2$. Indeed, on one hand, for any initial data $a(0)\in l^2$, there is a solution to (\ref{euler}) in $l^2$ satisfying the energy inequality
\begin{equation}\label{energy-euler}
\sum_{j=0}^\infty a_j^2(t)\leq \sum_{j=0}^\infty a_j^2(0) \ \ \ \ \forall \ \ t\geq 0.
\end{equation}
On the other hand, there exists a sequence $(A_j)_{j\geq 0}\in l^2$ such that the sequence 
\begin{equation}\notag
a_j(t)=\frac{A_j}{t_0-t}, \ \ \ \ t\in[0,t_0)
\end{equation}
is a self-similar solution to (\ref{euler}) in $l^2$. Note that $a_j$ are strictly increasing in time and hence this solution does not satisfy the energy inequality (\ref{energy-euler}). Nevertheless, the same authors showed in \cite{BFM-stoc2} that weak uniqueness can be recovered for the dyadic Euler model with some particular stochastic forcing.

We point out that there are non-unique weak solutions in $l^2$ for the deterministic model (\ref{sys-1}) as well since the Euler model (\ref{euler}) is a special case of (\ref{sys-1}). 
The purpose of this paper is to study the uniqueness problem of the dyadic MHD model with stochastic forcing. We expect to have weak uniqueness (uniqueness in law) with appropriate noise present in the model.

\medskip


Below we introduce various formulations of the dyadic MHD model (\ref{sys-1}) and its stochastic counter version. 

\subsection{Els\"asser form}
Let $P_j=a_j+b_j$ and $M_j=a_j-b_j$ for $j\geq 0$. Then $(P_j, M_j)$ satisfies
\begin{equation}\label{sys-pm}
\begin{split}
P_j'=&\ \lambda_{j-1}^\theta M_{j-1}P_{j-1}-\lambda_j^\theta M_jP_{j+1},\\
M_j'=&\ \lambda_{j-1}^\theta M_{j-1}P_{j-1}-\lambda_j^\theta P_jM_{j+1},
\end{split}
\end{equation}
for $j\geq 0$ and $P_{-1}=M_{-1}=0$. Compared to (\ref{sys-1}), (\ref{sys-pm}) contains lesson nonlinear terms, which is favorable while performing energy estimates in the analysis through out the text. 

\subsection{Stochastic system in Stratonovich form}
We consider the corresponding stochastic dyadic model in Stratonovich form with a special class of forcing
\begin{equation}\label{sys-St}
\begin{split}
d P_j=&\ (\lambda_{j-1}^\theta M_{j-1}P_{j-1}-\lambda_j^\theta M_jP_{j+1}) dt\\
&+\sigma \lambda_{j-1}^\theta P_{j-1} \circ dW_{j-1}^{p}-\sigma \lambda_j ^\theta P_{j+1}\circ d W_{j}^p,\\
d M_j=&\ (\lambda_{j-1}^\theta M_{j-1}P_{j-1}-\lambda_j^\theta P_jM_{j+1})dt\\
&+\sigma \lambda_{j-1}^\theta M_{j-1} \circ dW_{j-1}^m-\sigma \lambda_j^\theta M_{j+1}\circ d W_{j}^m
\end{split}
\end{equation}
for $j\geq 0$, $P_{-1}=M_{-1}=0$ and $\sigma\neq 0$. Here $W_j^p$ and $W_j^m$ are two independent families of Brownian motions. The motivation of the choice of such noise comes from the principle that the system formally conserves the energy. Indeed, one can check that 
\begin{equation}\notag
d\sum_{j=0}^\infty(P_j^2(t)+M_j^2(t))=0.
\end{equation}

\subsection{Stochastic system in It$\hat{\textbf{o}}$ form}
We also note the equivalent It$\hat{\textbf{o}}$ form of (\ref{sys-St}) is given by
\begin{equation}\label{sys-Ito}
\begin{split}
d P_j=&\ (\lambda_{j-1}^\theta M_{j-1}P_{j-1}-\lambda_j^\theta M_jP_{j+1}) dt\\
&+\sigma \lambda_{j-1}^\theta P_{j-1}  dW_{j-1}^{p}-\sigma \lambda_j^\theta P_{j+1} d W_{j}^p\\
&-\frac{\sigma^2}2(\lambda_j^{2\theta}+\lambda_{j-1}^{2\theta})P_jdt,\\
d M_j=&\ (\lambda_{j-1}^\theta M_{j-1}P_{j-1}-\lambda_j^\theta P_jM_{j+1})dt\\
&+\sigma \lambda_{j-1}^\theta M_{j-1}  dW_{j-1}^m-\sigma \lambda_j^\theta M_{j+1} d W_{j}^m\\
&-\frac{\sigma^2}2(\lambda_j^{2\theta}+\lambda_{j-1}^{2\theta})M_jdt,\\
\end{split}
\end{equation}
for $j\geq 1$, $P_{0}=M_{0}=0$ and $\sigma\neq 0$. 

\subsection{Stochastic system under Girsanov transform}

Denote 
\[U_j(t)=\frac{1}{\sigma}\int_0^t P_j(s)\, ds+W_j^m(t), \ \ V_j(t)=\frac{1}{\sigma}\int_0^t M_j(s)\, ds+W_j^p(t).\]
The It$\hat{\text{o}}$ form (\ref{sys-Ito}) can be written as
\begin{equation}\label{sys-G}
\begin{split}
d P_j=&\ \sigma\lambda_{j-1}^\theta P_{j-1}d V_{j-1}-\sigma\lambda_j^\theta P_{j+1}dV_j
-\frac{\sigma^2}2(\lambda_j^{2\theta}+\lambda_{j-1}^{2\theta})P_jdt,\\
d M_j=&\ \sigma\lambda_{j-1}^\theta M_{j-1}dU_{j-1}-\sigma\lambda_j^\theta M_{j+1}dU_j
-\frac{\sigma^2}2(\lambda_j^{2\theta}+\lambda_{j-1}^{2\theta})M_jdt,\\
\end{split}
\end{equation}
for $j\geq 1$, $P_{0}=M_{0}=0$ and $\sigma\neq 0$. Note that system (\ref{sys-G}) is formally linear.

\medskip

\subsection{Notion of solutions} 

Denote $\Omega=C([0,T]; \mathbb R)^{\mathbb N}$, $X(t)=(P(t), M(t))$, $X_j(t)=(P_j(t), M_j(t))$ and $Q=Q^p \times Q^m$.

\begin{Definition}\label{def-weak}
Given $X(0)=(p,m)\in l^2\times l^2$, a weak solution of (\ref{sys-Ito}) in $l^2\times l^2$ is a filtered probability space $(\Omega, F_t, Q)$, two families of independent Brownian motions $(W_j)_{j\geq0}$ on $(\Omega, F_t, Q)$, and an $l^2\times l^2$-valued stochastic process $(X_j)_{j\geq 1}$ on $(\Omega, F_t, Q)$ with continuous adapted components $P_j$ and $M_j$ such that
\begin{equation}\notag
\begin{split}
P_j(t)=&\ p_j+\int_0^t \left(\lambda_{j-1}^\theta M_{j-1}(s)P_{j-1}(s)-\lambda_j^\theta M_j(s)P_{j+1}(s)\right) \, ds\\
&+\int_0^t\sigma \lambda_{j-1}^\theta P_{j-1}(s)  dW^{p}_{j-1}(s)-\int_0^t\sigma \lambda_j^\theta P_{j+1}(s)\, d W^p_{j}(s)\\
&-\int_0^t\frac{\sigma^2}2(\lambda_j^{2\theta}+\lambda_{j-1}^{2\theta})P_j(s)\, ds,\\
M_j(t)=&\ m_j+\int_0^t \left(\lambda_{j-1}^\theta M_{j-1}(s)P_{j-1}(s)-\lambda_j^\theta P_j(s)M_{j+1}(s)\right)\,ds\\
&+\int_0^t\sigma \lambda_{j-1}^\theta M_{j-1}(s)  dW^m_{j-1}(s)-\int_0^t\sigma \lambda_j^\theta M_{j+1}(s)\, d W^m_{j}(s)\\
&-\int_0^t\frac{\sigma^2}2(\lambda_j^{2\theta}+\lambda_{j-1}^{2\theta})M_j(s)\, ds\\
\end{split}
\end{equation}
for $j\geq 1$ with $P_{0}=M_{0}=0$. This solution is denoted by $(\Omega, F_t, Q, W, X)$.
\end{Definition}

\begin{Definition}\label{def-exp}
A weak solution $(\Omega, F_t, Q, X)$ is said to be exponentially integrable if 
\begin{equation}\notag
\begin{split}
E^{Q^{p}}\left[\exp\left(\frac{1}{\sigma^2}\int_0^T\sum_{j=1}^{\infty}P_j^2(t)\, dt\right)\left(1+\int_0^TP_k^4(t)\, dt\right)^2\right]<\infty \ \ \forall \ \ k\geq 0,\\
E^{Q^{m}}\left[\exp\left(\frac{1}{\sigma^2}\int_0^T\sum_{j=1}^{\infty}M_j^2(t)\, dt\right)\left(1+\int_0^TM_k^4(t)\, dt\right)^2\right]<\infty \ \ \forall \ \ k\geq 0.
\end{split}
\end{equation}
\end{Definition}

\begin{Definition}\label{def-bounded}
A weak solution is said to be of class $L^\infty(\Omega\times [0,T], l^2\times l^2)$ if 
\begin{equation}\notag
\sum_{j=1}^\infty \left(P_j^2(t)+M_j^2(t)\right)\leq C \ \ \ \ \mbox{for a.e.} \ \ (\omega, t)\in \Omega\times [0,T]
\end{equation}
for a constant $C>0$. 
\end{Definition}

We note that a solution of class $L^\infty(\Omega\times [0,T], l^2\times l^2)$ is exponentially integrable. 

\begin{Definition}\label{def-unique}
System (\ref{sys-Ito}) is said to have weak uniqueness if it has uniqueness of the law of the process on the space $C([0,T]; \mathbb R)^{\mathbb N}$.
\end{Definition}

\begin{Definition}\label{def-weak2}
Given $X(0)=(p,m)\in l^2\times l^2$, a weak solution of (\ref{sys-G}) in $l^2\times l^2$ is a filtered probability space $(\Omega, F_t, \tilde Q)$, a sequence of independent Brownian motions $(\tilde W_j)_{j\geq1}=(V_j, U_j)_{j\geq1}$ on $(\Omega, F_t, \tilde Q)$, and an $l^2\times l^2$-valued stochastic process $(X_j)_{j\geq 1}$ on $(\Omega, F_t, \tilde Q)$ with continuous adapted components $P_j$ and $M_j$ such that
\begin{equation}\label{def-sol2}
\begin{split}
P_j(t)=&\ p_j+\int_0^t\sigma \lambda_{j-1}^\theta P_{j-1}(s)  dV_{j-1}(s)-\int_0^t\sigma \lambda_j^\theta P_{j+1}(s)\, d V_{j}(s)\\
&-\int_0^t\frac{\sigma^2}2(\lambda_j^{2\theta}+\lambda_{j-1}^{2\theta})P_j(s)\, ds,\\
M_j(t)=&\ m_j+\int_0^t\sigma \lambda_{j-1}^\theta M_{j-1}(s)  dU_{j-1}(s)-\int_0^t\sigma \lambda_j^\theta M_{j+1}(s)\, dU_{j}(s)\\
&-\int_0^t\frac{\sigma^2}2(\lambda_j^{2\theta}+\lambda_{j-1}^{2\theta})M_j(s)\, ds\\
\end{split}
\end{equation}
for $j\geq 1$ with $P_{0}=M_{0}=0$. This solution is denoted by $(\Omega, F_t, \tilde Q, \tilde W, X)$.
\end{Definition}

\medskip

\subsection{Main results} \label{sec-main}
For the choice of the particular stochastic forcing as in (\ref{sys-St}), we show uniqueness in law for the system.

\begin{Theorem}\label{thm-main}
For $p=(p_1, p_2, ...)\in l^2$ and $m=(m_1, m_2, ...)\in l^2$, there exists a weak $L^\infty(\Omega\times[0,T]; l^2\times l^2)$ solution to system (\ref{sys-St}) with initial data $X(0)=(p,m)$. Moreover, in the class of exponentially integrable solutions on $[0,T]$, there is weak uniqueness for system (\ref{sys-St}).
\end{Theorem}

The proof consists the four ingredients: (i) show the Stratonovich form (\ref{sys-St}) and the It$\hat{\textbf{o}}$ form (\ref{sys-Ito}) are equivalent; (ii) show the existence of weak $L^\infty$ solution to the It$\hat{\textbf{o}}$ form (\ref{sys-Ito}); (iii) in order to achieve the goal in step (ii), we apply Girsanov transform to (\ref{sys-Ito}) and study the resulted linear stochastic system (\ref{sys-G}); (iv) with existence and strong uniqueness results for the linear system (\ref{sys-G}), we establish the existence and weak uniqueness for the original nonlinear system (\ref{sys-Ito}). 

Although weak uniqueness is recovered, the weak solution obtained in Theorem \ref{thm-main} dissipates energy, which is defined by $\Upsilon(t) := \frac{1}{2}\sum_{n=1}^\infty \big( P^2_n(t)+M^2_n(t) \big)$.  This is given by the following theorem: 
\begin{Theorem} \label{anomalous_dissipation}
Let $\lambda_n = \lambda^n$ for some $\lambda> 1$. Given $p \in l^2$ and $m \in l^2$, let $X(t) = (P(t),M(t)) \in L^\infty(\Omega\times[0,T]; l^2\times l^2)$ be the unique weak solution obtained in Theorem \ref{thm-main}. Then, for any $t > 0$, 
\begin{equation*}
    Q \big ( \Upsilon(t) = \Upsilon(0) \big) < 1
\end{equation*}
and $\forall \epsilon > 0$, there exists $t > 0$ such that 
\begin{equation*}
    Q \big ( \Upsilon(t) < \epsilon \big) > 0
\end{equation*}
Furthermore, if the bound for the energy function $\Upsilon(t)$ is sufficiently small, then $\Upsilon(t)$ goes to 0 at least exponentially fast almost surely and in $L^1$.
\end{Theorem}

The key idea of proving Theorem \ref{anomalous_dissipation} is that we can recast the energy cascade mechanism of the system into a birth and death process, see more details in Section \ref{sec-energy}.

\bigskip

\section{Existence and uniqueness of the auxiliary linear system}\label{sec-linear}

This section is devoted to the the well-posedness of the linear system (\ref{sys-G}). 
\begin{Theorem}\label{thm-G}
Let $p\in l^2$ and $m\in l^2$. For system (\ref{sys-G}) with initial data $X(0)=(p,m)$, there exists a unique solution in $L^\infty(\Omega\times [0,T]; l^2\times l^2)$ with continuous components. 
\end{Theorem}

\pf
We first show the existence of a solution by applying Galerkin's approximating method. For each integer $N\geq 1$, we consider the finite dimensional stochastic system
\begin{equation}\label{sys-Gcut}
\begin{split}
d P_j^{(N)}=&\ \sigma\lambda_{j-1}^\theta P_{j-1}^{(N)}d V_{j-1}-\sigma\lambda_j^\theta P_{j+1}^{(N)}dV_j
-\frac{\sigma^2}2(\lambda_j^2+\lambda_{j-1}^2)P_j^{(N)}dt,\\
d M_j^{(N)}=&\ \sigma\lambda_{j-1}^\theta M_{j-1}^{(N)}dU_{j-1}-\sigma\lambda_j^\theta M_{j+1}^{(N)}dU_j
-\frac{\sigma^2}2(\lambda_j^2+\lambda_{j-1}^2)M_j^{(N)}dt,\\
P_j^{(N)}(0)=&\ p_j, \ \ \ M_j^{(N)}(0)= m_j, 
\end{split}
\end{equation}
with 
\[P_0^{(N)}(t)=P_{N+1}^{(N)}(t)=M_0^{(N)}(t)=M_{N+1}^{(N)}(t)=0.\]
By standard theory of ODE with stochastic force, there exists a unique strong solution on $[0,\infty)$ to system (\ref{sys-Gcut}). It$\hat{\text{o}}$'s formula gives
\begin{equation}\label{ito-pm}
\begin{split}
\frac12 d\left(P_j^{(N)}\right)^2=&\ P_j^{(N)}d P_j^{(N)}+\frac12d\left[P_j^{(N)}, P_j^{(N)}\right]_t,\\
\frac12 d\left(M_j^{(N)}\right)^2=&\ M_j^{(N)}d M_j^{(N)}+\frac12d\left[M_j^{(N)}, M_j^{(N)}\right]_t. 
\end{split}
\end{equation}
On the other hand, it follows from system (\ref{sys-Gcut}) that
\begin{equation}\label{covar}
\begin{split}
\left[P_j^{(N)}, P_j^{(N)}\right]_t=\left[P_j^{(N)}\right]_t=&\int_0^t\sigma^2\lambda_{j-1}^{2\theta}\left(P_{j-1}^{(N)}(s)\right)^2
+\sigma^2\lambda_{j}^{2\theta}\left(P_{j+1}^{(N)}(s)\right)^2\, ds\\
\left[M_j^{(N)}, M_j^{(N)}\right]_t=\left[M_j^{(N)}\right]_t=&\int_0^t\sigma^2\lambda_{j-1}^{2\theta}\left(M_{j-1}^{(N)}(s)\right)^2
+\sigma^2\lambda_{j}^{2\theta}\left(M_{j+1}^{(N)}(s)\right)^2\, ds.
\end{split}
\end{equation}
Combining (\ref{ito-pm}) and (\ref{covar}) leads to
\begin{equation}\notag
\begin{split}
\frac12 d\left(P_j^{(N)}\right)^2=&\ P_j^{(N)}d P_j^{(N)}+\frac12\sigma^2\lambda_{j-1}^{2\theta}\left(P_{j-1}^{(N)}\right)^2dt+\frac12\sigma^2\lambda_{j}^{2\theta}\left(P_{j+1}^{(N)}\right)^2dt,\\
\frac12 d\left(M_j^{(N)}\right)^2=&\ M_j^{(N)}d M_j^{(N)}+\frac12\sigma^2\lambda_{j-1}^{2\theta}\left(M_{j-1}^{(N)}\right)^2dt+\frac12\sigma^2\lambda_{j}^{2\theta}\left(M_{j+1}^{(N)}\right)^2dt. 
\end{split}
\end{equation}
Adding the last two equations and taking sum for $0\leq j\leq N$ we obtain
\begin{equation}\notag
\begin{split}
\frac12 d\sum_{j=0}^N\left(P_j^{(N)}\right)^2
=&\sum_{j=0}^N \sigma\lambda_{j-1}^{2\theta} P_{j-1}^{(N)}P_j^{(N)}dV_{j-1}-\sum_{j=0}^N \sigma\lambda_{j}^{2\theta} P_{j}^{(N)}P_{j+1}^{(N)}dV_{j}\\
&-\sum_{j=0}^N\frac{\sigma^2}2(\lambda_j^2+\lambda_{j-1}^2)\left(P_j^{(N)}\right)^2dt\\
&+\sum_{j=0}^N\frac12\sigma^2\lambda_{j-1}^{2\theta}\left(P_{j-1}^{(N)}\right)^2dt+\sum_{j=0}^N\frac12\sigma^2\lambda_{j}^{2\theta}\left(P_{j+1}^{(N)}\right)^2dt\\
=&-\frac12\sigma^2\lambda_{N}^{2\theta}\left(P_{N}^{(N)}\right)^2dt,
\end{split}
\end{equation}
and similarly
\begin{equation}\notag
\begin{split}
\frac12 d\sum_{j=0}^N\left(M_j^{(N)}\right)^2
=&-\frac12\sigma^2\lambda_{N}^{2\theta}\left(M_{N}^{(N)}\right)^2dt.
\end{split}
\end{equation}
We note that 
\begin{equation}\notag
-\frac12\sigma^2\lambda_{N}^{2\theta}\left(P_{N}^{(N)}\right)^2<0, \ \ -\frac12\sigma^2\lambda_{N}^{2\theta}\left(M_{N}^{(N)}\right)^2<0.
\end{equation}
Therefore know
\begin{equation}\notag
\sum_{j=0}^N\left(P_j^{(N)}(t)\right)^2\leq \sum_{j=0}^N p_j^2, \ \ \sum_{j=0}^N\left(M_j^{(N)}(t)\right)^2\leq \sum_{j=0}^N m_j^2, \ \ \tilde{Q}-a.s.
\end{equation}
As a consequence, we know that there exists a subsequence of $\left(P_j^{(N)}, M_j^{(N)}\right)_{j\geq 1}$ which converges to certain 
$\left(P_j, M_j\right)_{j\geq 1}\in L^p(\Omega \times [0,T]; l^2\times l^2))$; and the convergence is weak star in $L^\infty(\Omega \times [0,T]; l^2\times l^2)$. Following standard arguments for equations with monotone operators (see \cite{PM}), the subspace of progressively measurable processes of $L^p(\Omega \times [0,T]; l^2\times l^2)$ is strongly closed and hence $\left(P_j, M_j\right)_{j\geq 1}$ is progressively measurable. Furthermore, the stochastic integrals in (\ref{def-sol2}) are strongly continuous linear operators from the subspace of progressively measurable processes of $L^2(\Omega \times [0,T]; l^2\times l^2)$ to $L^2(\Omega)$. Therefore, integrating (\ref{sys-Gcut}) and passing to the limit in each integral gives us the equations of (\ref{def-sol2}) satisfied by the limit $\left(P_j, M_j\right)_{j\geq 1}$. In view of the integral equations of (\ref{def-sol2}), we know that there is a version of $\left(P_j, M_j\right)_{j\geq 1}$ such that all components are continuous. This completes the proof of the existence. The uniqueness is an immediate consequence of Proposition \ref{prop} below, since any process in $L^\infty(\Omega \times [0,T]; l^2\times l^2)$ belongs to the class $\mathcal K$ defined in (\ref{K}).

\cbdu

\medskip

Denote the class of random variables 
\begin{equation}\label{K}
\begin{split}
\mathcal K=& \left\{ (P_j, M_j)_{j\geq 1}:  \int_0^T E^{Q^p}[P_j^4(t)]\, dt + \int_0^T E^{Q^m}[M_j^4(t)]\, dt<\infty \ \ \forall \ \ j\geq0, \right.\\
&  \ \ \ \ \ \ \ \ \ \ \ \ \ \ \ \ \ \ \ \ \ \ \mbox{and} \ \ \left.  \lim_{j\to\infty}\int_0^T E^{Q^p}[P_j^2(t)]\, dt + \int_0^T E^{Q^m}[M_j^2(t)]\, dt = 0.
\right\}
\end{split}
\end{equation}

\begin{Proposition}\label{prop}
Associated with initial data $(p,m)\in l^2\times l^2$, there is at most one solution to system (\ref{sys-G}) in the class $\mathcal K$ on $[0,T]$.
\end{Proposition}
\pf 
It is sufficient to show that the only solution in class $\mathcal K$ of (\ref{sys-G}) with zero initial condition is the null process. Thus with $p=m=0$ we have from (\ref{def-sol2}) by It$\hat{\text{o}}$ formula  
\begin{equation}\label{est-energy1}
\begin{split}
\frac12 d P_j^2=&\ P_j dP_j+\frac12 d[P_j]_t\\
=&-\frac12\sigma^2(\lambda_j^{2\theta}+\lambda_{j-1}^{2\theta})P_j^2 dt+\frac12\sigma^2(\lambda_{j-1}^{2\theta}P_{j-1}^2+\lambda_j^{2\theta}P_{j+1}^2) dt+ d \tilde V_j,\\
\frac12 d M_j^2=&\ M_j dM_j+\frac12 d[M_j]_t\\
=&-\frac12\sigma^2(\lambda_j^{2\theta}+\lambda_{j-1}^{2\theta})M_j^2 dt+\frac12\sigma^2(\lambda_{j-1}^{2\theta}M_{j-1}^2+\lambda_j^{2\theta}M_{j+1}^2) dt+ d \tilde U_j
\end{split}
\end{equation}
with 
\begin{equation}\notag
\begin{split}
\tilde V_j(t)=&\int_0^t \sigma \lambda_{j-1}^\theta P_{j-1}(s)P_j(s) dV_{j-1}(s)-\int_0^t \sigma \lambda_{j}^\theta P_{j}(s)P_{j+1}(s) dV_{j}(s),\\
\tilde U_j(t)=&\int_0^t \sigma \lambda_{j-1}^\theta M_{j-1}(s)M_j(s) dU_{j-1}(s)-\int_0^t \sigma \lambda_{j}^\theta M_{j}(s)M_{j+1}(s) dU_{j}(s).
\end{split}
\end{equation}
The first condition of class $\mathcal K$ in (\ref{K}) implies that $E^{Q^p}[P_j^2(t)]$ and $E^{Q^m}[M_j^2(t)]$ are finite and continuous in $t$.
It also implies $\tilde V_j$ and $\tilde U_j$ are martingales for all $j\geq 1$ and hence 
\[E^{Q^p}[\tilde V_j]=0, \ \ \ E^{Q^m}[\tilde U_j]=0 \ \ \ \forall \ \ j\geq 1.\]
Therefore it follows from (\ref{est-energy1}) that for all $j\geq 1$
\begin{equation}\label{est-energy2}
\begin{split}
E^{Q^p}[P_j^2(t)]=&-\sigma^2(\lambda_j^{2\theta}+\lambda_{j-1}^{2\theta})\int_0^t E^{Q^p}[P_j^2(s)]ds\\
&+\sigma^2\lambda_{j-1}^{2\theta} \int_0^t E^{Q^p}[P_{j-1}^2(s)]ds+\sigma^2\lambda_{j}^{2\theta} \int_0^t E^{Q^p}[P_{j+1}^2(s)]ds,\\
E^{Q^m}[M_j^2(t)]=&-\sigma^2(\lambda_j^{2\theta}+\lambda_{j-1}^{2\theta})\int_0^t E^{Q^m}[M_j^2(s)]ds\\
&+\sigma^2\lambda_{j-1}^{2\theta} \int_0^t E^{Q^m}[M_{j-1}^2(s)]ds+\sigma^2\lambda_{j}^{2\theta} \int_0^t E^{Q^m}[M_{j+1}^2(s)]ds.
\end{split}
\end{equation}
Further we deduce from (\ref{est-energy2}) that
\begin{equation}\notag
\begin{split}
\int_0^t E^{Q^p}[(P_{j+1}^2(s)-P_j^2(s))]ds \geq \frac{\lambda_{j-1}^{2\theta}}{\lambda_j^{2\theta}}\int_0^t E^{Q^p}[(P_{j}^2(s)-P_{j-1}^2(s))]ds,\\
\int_0^t E^{Q^m}[(M_{j+1}^2(s)-M_j^2(s))]ds \geq \frac{\lambda_{j-1}^{2\theta}}{\lambda_j^{2\theta}}\int_0^t E^{Q^m}[(M_{j}^2(s)-M_{j-1}^2(s))]ds.
\end{split}
\end{equation}
In view of $P_0\equiv0$ and $M_0\equiv0$, we have by iteration that
\begin{equation}\notag
\int_0^t E^{Q^p}[(P_{j+1}^2(s)-P_j^2(s))]ds\geq 0, \ \ \int_0^t E^{Q^m}[(M_{j+1}^2(s)-M_j^2(s))]ds\geq 0 \ \ \ \forall \ \ j\geq 1.
\end{equation}
It thus follows from the second condition of class $\mathcal K$ of (\ref{K}) that 
\begin{equation}\notag
\int_0^t E^{Q^p}[P_j^2(s)]ds= 0, \ \ \int_0^t E^{Q^m}[M_j^2(s)]ds= 0 \ \ \ \forall \ \ j\geq 1.
\end{equation}
Therefore the uniqueness in class $\mathcal K$ is obtained.
\cbdu
\bigskip

\section{Proof of the Theorem \ref{thm-main}}\label{sec-proof}

We present a proof for Theorem \ref{thm-main} in this section. As outlined in Subsection \ref{sec-main}, we first show that a weak solution of the It$\hat{\text{o}}$ form (\ref{sys-Ito}) is also a weak solution of the Stratonovich form (\ref{sys-St}). Namely, we will prove the lemma below.

\begin{Lemma}\label{le-equ}
Let $X(t)=(P(t), M(t))$ be a weak solution of system (\ref{sys-Ito}). Then $X$ satisfies the Stratonovich system (\ref{sys-St}). 
\end{Lemma}
\pf
Note that the processes $(P_j)_{j\geq 1}$ and $(M_j)_{j\geq 1}$ are continuous semimartingales for all $j\geq 1$. Thus the Stratonovich integrals arising in (\ref{sys-St}) are well-defined. Moreover, we have
\begin{equation}\notag
\int_0^t P_{j-1}(s)\circ d W_{j-1}(s) =\int_0^t P_{j-1}(s)d W_{j-1}(s)+\frac12[P_{j-1}, W_{j-1}]_t.
\end{equation}
It follows from the equation of $P_{j-1}$ 
\begin{equation}\notag
\begin{split}
d{P_{j-1}} =&\ (\lambda_{j-2}^\theta M_{j-2}P_{j-2}-\lambda_{j-1}^\theta M_{j-1}P_j) dt\\
&+\sigma \lambda_{j-2}^\theta P_{j-2}\circ d {W_{j-2}}-\sigma \lambda_{j-1}^\theta P_j\circ dW_{j-1}
\end{split}
\end{equation}
and the independence of the Brownian motions that the joint quadratic variation of $P_{j-1}$ and $W_{j-1}$ is given by
\begin{equation}\notag
[P_{j-1}, W_{j-1}]_t=-\int_0^t\sigma \lambda_{j-1}^\theta P_j(s) ds. 
\end{equation}
Therefore we obtain
\begin{equation}\notag
\int_0^t\sigma \lambda_{j-1}^\theta P_{j-1}(s)\circ d W_{j-1}(s) = \int_0^t\sigma \lambda_{j-1}^\theta P_{j-1}(s) d W_{j-1}(s)-\frac12\int_0^t\sigma^2\lambda_{j-1}^{2\theta} P_j(s) ds.
\end{equation}
Analogously, we have
\begin{equation}\notag
\begin{split}
\int_0^t\sigma \lambda_{j}^\theta P_{j+1}(s)\circ d W_{j}(s) = &\int_0^t\sigma \lambda_{j}^\theta P_{j+1}(s) d W_{j}(s)+\frac12\int_0^t\sigma^2\lambda_{j}^{2\theta} P_j(s) ds, \\
\int_0^t\sigma \lambda_{j-1}^\theta M_{j-1}(s)\circ d W_{j-1}(s) =& \int_0^t\sigma \lambda_{j-1}^\theta M_{j-1}(s) d W_{j-1}(s)-\frac12\int_0^t\sigma^2\lambda_{j-1}^{2\theta} M_j(s) ds,\\
\int_0^t\sigma \lambda_{j}^\theta M_{j+1}(s)\circ d W_{j}(s) =& \int_0^t\sigma \lambda_{j}^\theta M_{j+1}(s) d W_{j}(s)+\frac12\int_0^t\sigma^2\lambda_{j}^{2\theta} M_j(s) ds. \\
\end{split}
\end{equation}
Therefore we conclude that a weak solution $X(t)=(P(t), M(t))$ of (\ref{sys-Ito}) satisfies (\ref{sys-St}) as well.

\cbdu

\medskip

Next we show the existence of weak $L^\infty$ solution to the It$\hat{\text{o}}$ system (\ref{sys-Ito}) by applying results established in Section \ref{sec-linear} for the linearized system (\ref{sys-G}). 

\begin{Lemma}\label{le-ex-Ito}
For $p=(p_1, p_2, ...)\in l^2$ and $m=(m_1, m_2, ...)\in l^2$, there exists a weak $L^\infty(\Omega\times[0,T]; l^2\times l^2)$ solution to system (\ref{sys-Ito}) associated with initial data $X(0)=(p,m)$. 
\end{Lemma}
\pf
Let $(\Omega, F_t, \tilde Q, \tilde W, X)$ be the unique solution in $L^\infty(\Omega \times [0,T]; l^2\times l^2)$ to the linear system (\ref{sys-G}) associated with initial data $(p,m)$ obtained in Theorem \ref{thm-G}. Note that the solution is exponentially integrable and hence
\begin{equation}\notag
E^Q\left[\int_0^T \sum_{j=1}^\infty (P_j^2(s)+M_j^2(s)) ds< \infty \right].
\end{equation}
Denote the processes 
\begin{equation}\notag
L_{1,t}=-\frac1\sigma \sum_{j=1}^\infty\int_0^t P_j(s) d W_j(s), \ \ \ L_{2,t}=-\frac1\sigma \sum_{j=1}^\infty\int_0^t M_j(s) d W_j(s)
\end{equation}
which are well-defined and are martingales. Moreover, we have
\begin{equation}\notag
[L_1, L_1]_t=\frac{1}{\sigma^2}\int_0^t \sum_{j=1}^{\infty} P_j^2(s) ds, \ \ [L_2, L_2]_t=\frac{1}{\sigma^2}\int_0^t \sum_{j=1}^{\infty} M_j^2(s) ds.
\end{equation}
Similarly, the processes
\begin{equation}\notag
Z_{1,t}=\frac1\sigma \sum_{j=1}^\infty\int_0^t P_j(s) d U_j(s), \ \ \ Z_{2,t}=\frac1\sigma \sum_{j=1}^\infty\int_0^t M_j(s) d V_j(s)
\end{equation}
with
\[U_j(t)=\frac{1}{\sigma}\int_0^t P_j(s)\, ds+W_j(t), \ \ V_j(t)=\frac{1}{\sigma}\int_0^t M_j(s)\, ds+W_j(t)\]
are well-defined and are martingales, and the quadratic variations are 
\begin{equation}\notag
[Z_1, Z_1]_t=\frac{1}{\sigma^2}\int_0^t \sum_{j=1}^{\infty} P_j^2(s) ds, \ \ [Z_2, Z_2]_t=\frac{1}{\sigma^2}\int_0^t \sum_{j=1}^{\infty} M_j^2(s) ds.
\end{equation}
Recall $\tilde Q=(\tilde Q_1, \tilde Q_2)$ and $\tilde W=(U, V)$. The probability measures $\tilde Q_1$ and $\tilde Q_2$ are associated with $U$ and $V$ respectively.  We introduce the new measure $Q$ on $(\Omega, F_t)$ satisfying 
\begin{equation}\notag
\begin{split}
\frac{dQ^p}{d\tilde Q_1}=&\ \exp \left(Z_{1,T}-\frac12[Z_1, Z_1]_T\right),\\
\frac{dQ^m}{d\tilde Q_2}=&\ \exp \left(Z_{2,T}-\frac12[Z_2, Z_2]_T\right).\\
\end{split}
\end{equation}
In view of Novikov's criterion and the fact
\begin{equation}\notag
E\left[\exp\left(\frac{1}{2\sigma^2}\int_0^T\sum_{j=1}^\infty P_j^2(t)dt\right)\right]<\infty, \ \ \ E\left[\exp\left(\frac{1}{2\sigma^2}\int_0^T\sum_{j=1}^\infty M_j^2(t)dt\right)\right]<\infty,
\end{equation}
we know that $\exp \left(Z_{1,T}-\frac12[Z_1, Z_1]_T\right)$ and $\exp \left(Z_{2,T}-\frac12[Z_2, Z_2]_T\right)$ are strictly positive martingales. Therefore the measures $Q$, $\tilde Q_1$ and $\tilde Q_2$ are equivalent on $F_T$. Hence we have
\begin{equation}\label{measure-Q12}
\begin{split}
\frac{d\tilde Q_1}{dQ^p}=&\ \exp \left(L_{1,T}-\frac12[L_1, L_1]_T\right)=\exp \left(-Z_{1,T}+\frac12[Z_1, Z_1]_T\right),\\
\frac{d\tilde Q_2}{dQ^m}=&\ \exp \left(L_{2,T}-\frac12[L_2, L_2]_T\right)=\exp \left(-Z_{2,T}+\frac12[Z_2, Z_2]_T\right).\\
\end{split}
\end{equation}
The precesses under $Q$ 
\begin{equation}\notag
W_j(t)=U_j(t)-\frac{1}{\sigma}\int_0^t P_j(s)ds=V_j(t)-\frac{1}{\sigma}\int_0^t M_j(s)ds
\end{equation}
are a sequence of independent Brownian motions. Therefore $(\Omega, F_t, Q, W, X)$ is a solution to the nonlinear system (\ref{sys-Ito}). Moreover, the solution is in $L^\infty(\Omega\times[0,T]; l^2\times l^2)$ since the measures $Q$, $\tilde Q_1$ and $\tilde Q_2$ are equivalent on $F_T$.

\cbdu

\medskip

\begin{Lemma}\label{le-existence3}
Let $(\Omega, F_t, Q, W, X)$ be an exponentially integrable solution of the nonlinear system (\ref{sys-Ito}). Then it is a solution to the linear system (\ref{sys-G}) with processes 
\begin{equation}\notag
U_j(t)=\frac{1}{\sigma}\int_0^t P_j(s)\, ds+W_j(t), \ \ V_j(t)=\frac{1}{\sigma}\int_0^t M_j(s)\, ds+W_j(t)
\end{equation}
being independent Brownian motions on $(\Omega, F_T, \tilde Q_1)$ and $(\Omega, F_T, \tilde Q_2)$ respectively. The measures $\tilde Q_1$ and $\tilde Q_2$ are defined by (\ref{measure-Q12}). Moreover, the process $X$ on 
\[(\Omega, F_T, \tilde Q, \tilde W):=(\Omega, F_T, (\tilde Q_1, \tilde Q_2), (U, V))\] is in class $\mathcal K$ as defined in (\ref{K}).
\end{Lemma}

\pf
The proof of the existence part can be obtained analogously as in the proof of Lemma \ref{le-ex-Ito}. We only need to verify that the process $X$ is in class $\mathcal K$. 
\begin{equation}\notag
\begin{split}
&E^{\tilde Q_1} \left[\int_0^T P_j^4(t)dt\right]\\
=&\ E^Q \left[\exp\left(L_{1,T}-\frac12[L_1, L_1]_T\right)\int_0^T P_j^4(t)dt\right]\\
\leq &\ E^Q\left[\exp(2L_{1,T}-2[L_1, L_1]_T)\right]^{\frac12} E^Q\left[\left(\int_0^T P_j^4(t)dt\right)^2\exp([L_1, L_1]_T)\right]^{\frac12}.
\end{split}
\end{equation}
Note that
\begin{equation}
E^Q\left[\exp(2L_{1,T}-2[L_1, L_1]_T)\right]=1
\end{equation}
by applying Girsanov's theorem to the martingale $2L_{1,t}$. 

\cbdu

\medskip

{\textbf {Finishing the proof of Theorem \ref{thm-main}: }}  We are left to show weak uniqueness. Following standard arguments, we assume that $(\Omega^{(1)}, F_t^{(1)}, Q^{(1)}, W^{(1)}, X^{(1)})$ and $(\Omega^{(2)}, F_t^{(2)}, Q^{(2)}, W^{(2)}, X^{(2)})$ are two exponentially integrable solutions of system (\ref{sys-Ito}) with the same initial data $X(0)=(p,m)\in l^2\times l^2$. Then for $i = 1,2$, we have 
\\
\begin{equation}\label{sys-uniqueness}
\begin{split}
   d P^{(i)}_j=&\ \sigma\lambda_{j-1}^\theta P^{(i)}_{j-1}d V^{(i)}_{j-1}-\sigma\lambda_j^\theta P^{(i)}_{j+1}dV^{(i)}_j
-\frac{\sigma^2}2(\lambda_j^{2\theta}+\lambda_{j-1}^{2\theta})P^{(i)}_jdt,\\
d M^{(i)}_j=&\ \sigma\lambda_{j-1}^\theta M^{(i)}_{j-1}dU^{(i)}_{j-1}-\sigma\lambda_j^\theta M^{(i)}_{j+1}dU^{(i)}_j
-\frac{\sigma^2}2(\lambda_j^{2\theta}+\lambda_{j-1}^{2\theta})M^{(i)}_jdt,\\
\end{split}
\end{equation}
where,
      \[U^{(i)}_j(t)=\frac{1}{\sigma}\int_0^t P^{(i)}_j(s)\, ds+W_j^{(i),m(t)}, \ \ V_j(t)=\frac{1}{\sigma}\int_0^t M_j(s)\, ds+W_j^{(i),p(t)}.\]
are sequences of Brownian motions on $(\Omega^{(i)}, F_t^{(i)}, \tilde{Q}^{(i)})$ in which the measure $\tilde{Q}^{(i)}$ are given by $(\ref{measure-Q12})$. Since in Theorem $(\ref{thm-G})$ we have shown that there exists a unique solution in $L^\infty(\Omega \times [0,T];l^2\times l^2)$ in the system $(\ref{sys-G})$, by Yamata-Watanabe theorem, the system has uniquness in law of the processs on $C([0,T];\mathbb{R})^\mathbb{N}$.
Then given $n\in \mathbb{N}$, $t_i \in [0,T]$, $i = 1,...,n$, and $f:(l^2\times l^2)^n \to \mathbb{R}$ a measurable bounded function, by (\ref{measure-Q12}) we have that 

\begin{equation}\notag
    \begin{split}
       & E^{Q^{(i)}}\Big[f\big(P^{(i)}(t_1),...,P^{(i)}(t_n)\big) \Big] \\
        =  &E^{\tilde{Q}_1^{(i)}}\Big[\exp \left(Z^{(i)}_{1,T}-\frac12[Z^{(i)}_1, Z^{(i)}_1]_T\right)f\big(P^{(i)}(t_1),...,P^{(i)}(t_n)\big) \Big]\\
    \end{split}
\end{equation}
where 
\begin{equation}\notag
    Z^{(i)}_{1,t} =  \frac1\sigma \sum_{j=1}^\infty\int_0^t P_j(s) d U_j(s)
\end{equation}
Therefore,
\begin{equation}\notag
    \begin{split}
       & E^{Q^{(1)}}\Big[f\big(P^{(1)}(t_1),...,P^{(1)}(t_n)\big) \Big] \\
        =  &E^{\tilde{Q}_1^{(1)}}\Big[\exp \left(Z^{(1)}_{1,T}-\frac12[Z^{(1)}_1, Z^{(1)}_1]_T\right)f\big(P^{(1)}(t_1),...,P^{(1)}(t_n)\big) \Big]\\
        =  &E^{\tilde{Q}_1^{(2)}}\Big[\exp \left(Z^{(2)}_{1,T}-\frac12[Z^{(2)}_1, Z^{(2)}_1]_T\right)f\big(P^{(2)}(t_1),...,P^{(2)}(t_n)\big) \Big]\\
        = &E^{Q^{(2)}}\Big[f\big(P^{(2)}(t_1),...,P^{(2)}(t_n)\big) \Big] 
    \end{split}
\end{equation}\\
The third equality above can be deduced by enlarging the systems $(\ref{sys-uniqueness})$  with 
\begin{equation}\notag
\begin{split}
    dZ^{(i)}_{1,t} &=  \frac1\sigma \sum_{j=1}^\infty P^{(i)}_j(s) d U^{(i)}_j(s)\\
    dZ^{(i)}_{2,t} &=  \frac1\sigma \sum_{j=1}^\infty M^{(i)}_j(s) d V^{(i)}_j(s)
\end{split}
\end{equation}
The enlarged system also has strong uniqueness, so again by Yamata-Watanabe Theorem, it has weak uniqueness. Repeating the above procedure on the second component $M^{(i)}$ we derive the weak uniqueness.  
\bigskip

\section{Anomalous dissipation of energy}
\label{sec-energy}
\vspace{0.5cm}

The proof of Theorem \ref{anomalous_dissipation} follows mainly from techniques about birth and death process, similarly as for the Euler dyadic model in \cite{BFM-stoc1}. 

\subsection{Birth and Death process}
We first introduce an associated birth and death process. Let \\
\begin{equation}\notag
\begin{split}
     e_j(t) &= \frac{1}{||x||^2}\Big(E^{\tilde{Q}}[P^2_j(t)+M^2_j(t)] \Big)\\
     ||x||^2_{2} &= ||p||^2_{2}+||m||^2_{2}
\end{split}
\end{equation}\\
and we denote $e(t) := (e_j(t))_{j \geq 1} $, $t \geq 0 $, $e_0(t)\equiv 0$. In addition, we let \\
\begin{equation}\notag
    \mu_j = \sigma^2\lambda^{2\theta}_{j-1} \hspace{1cm}   \nu_j = \sigma^2\lambda^{2\theta}_{j}
\end{equation}
From $(\ref{est-energy2})$ we derive that \\
\begin{equation}\notag
    \begin{split}
        \frac{d}{dt}E^{\tilde{Q}}[P_j^2(t)]=&-\sigma^2(\lambda_j^{2\theta}+\lambda_{j-1}^{2\theta}) E^{Q}[P_j^2(t)]\\
&+\sigma^2\lambda_{j-1}^{2\theta}  E^{\tilde{Q}}[P_{j-1}^2(t)]+\sigma^2\lambda_{j}^{2\theta}  E^{Q}[P_{j+1}^2(t)],\\
\frac{d}{dt}E^{\tilde{Q}}[M_j^2(t)]=&-\sigma^2(\lambda_j^{2\theta}+\lambda_{j-1}^{2\theta}) E^{\tilde{Q}}[M_j^2(t)]\\
&+\sigma^2\lambda_{j-1}^{2\theta}E^{\tilde{Q}}[M_{j-1}^2(t)]
+\sigma^2\lambda_{j}^{2\theta}  E^{\tilde{Q}}[M_{j+1}^2(t)].  
    \end{split}
\end{equation}\\
sum up the above two equations and rewrite it in terms of $e(t)$ \\ 
\begin{equation}\label{birth_death_process}
    \begin{split}
        \frac{d}{dt}e_j(t) &= -(\nu_j+\mu_{j}) e_j(t)
+\nu_{j-1} e_{j-1}(t)+\mu_{j+1}e_{j+1}(t)\\
        e_j(0) &= \frac{p^2_j+m^2_j}{||x||^2_{2}}
    \end{split}
\end{equation}\\
The system $(\ref{birth_death_process})$ can be written as $\frac{d}{dt}e(t) = e(t)$A,
where $A$ is an infinite matrix such that $A = \{A_{m,n}\}$, $A_{m,m} = -(\nu_m+\mu_m)$, $A_{m,m+1} = \mu_{m+1}$, $A_{m,m-1} = \nu_{m-1}$ and $A_{m,n} = 0$ otherwise. We call this matrix $A$ a q-matrix, and the processes give by such q-matrix $A$ $\textit{birth and death processes}$. Our goal is to find a (unique) process $\eta_t$ on a new probability space $(Y,\mathcal{Y},\mathbb{P})$ such that $e_j(t) = \mathbb{P}(\eta_t = j)$. The construction is as followed: Given a probability space $(Y,\mathcal{Y},\mathbb{P})$, let $\eta_t$ be a continuous time Markov chain on positive integers, with initial distribution
\begin{equation}\notag
\mathbb{P}(\eta_0 = j) = e_j(0), \ \ \ \ \  j = 1,2,3...
\end{equation}
    
 The jump rates of $\eta_t$ are given by the entries of $A$, i.e. $\eta_t$ waits in a state $j$ for an exponential time with rate $\mu_j+\nu_j$ before going to $j+1$ or $j-1$ with the probability $\Pi_{j,j+1}$ and $\Pi_{j,j-1}$ respectively, where 
 \begin{equation}\notag
     \Pi_{j,j+1} = \frac{\nu_j}{\mu_j+\nu_j}, \ \ \ \ \   \Pi_{j,j-1} = \frac{\mu_j}{\mu_j+\nu_j}
 \end{equation}
 
 Since $A$ is stable , i.e. $-A_{m,m} \leq \infty$ for all $m$, by Theorem 2.1 of Chapter 2 in \cite{And}, there exists a process whose law is a solution to both forward equation ($\frac{d}{dt}e(t) = e(t)$A) and backward equations ($\frac{d}{dt}e(t) = Ae(t)$). Such process is called a $\textit{minimal solution}$. \\\\
 Denote $\tau \in [0,\infty]$ the first time such that in $[0,\tau)$ the process experience infinitely many jumps. If the minimal solution of a q-matrix is honest, which means $\mathbb{P}(\tau < \infty) = 0$, then the minimal solutions for forward and backward equations both have uniqueness. (Theorem 2.2.2 in \cite{And}) And such q-matrix is called $\textit{regular}$ (see definition in \cite{And}, page 81 for more detail). Unfortunately, it turns out that the minimal solution constructed above is dishonest. Nevertheless, the solution to the forward equations is still unique, even though the backward equations have infinitely many solutions. We will state the key theorems below, from which the uniqueness results follow. 

\begin{Theorem}(Theorem 3.2.2 in \cite{And})\label{backward_equ_uni_thm}
Let Q be the q-matrix of the birth and death process associated with (\ref{birth_death_process}), and assume that $\nu_n >0$ for all $n \geq 1$. Define
\begin{equation}\label{quantityR}
    R = \sum_{k=1}^{\infty} \Big(\frac{1}{\nu_k} +\frac{\mu_k}{\nu_k\nu_{k-1}}+\frac{\mu_k\mu_{k-1}}{\nu_k\nu_{k-1}\nu_{k-2}}+...+\frac{\mu_k\mu_{k-1}...\mu_2}{\nu_k\nu_{k-1}\nu_{k-2}...\nu_1}\Big)
\end{equation}
Then the equation $Qx=\lambda x,\: 0 \leq x \leq 1 \: (0 \leq x_i \leq 1, i=0,1,2,... )$, has only the trivial solution x = 0 if and only if $R = \infty$. This implies that the minimal solution is the unique solution of the backward equation if and only if $R = \infty$. 
\end{Theorem}\vspace{0.5cm}
\begin{Corollary}
The q-matrix A is is not regular if and only if $\sum_{k} k\lambda^{-2\theta}_k < \infty$. Therefore the minimal solution is dishonest.
\end{Corollary}
\pf
It follows from the Corollary 2.2.5 of \cite{And} that, since $A$ is conservative (sum of each row of A is 0), the minimal solution is unique if and only if it is honest. By Theorem \ref{backward_equ_uni_thm}, the minimal solution of the backward equation is unique if and only if the quantity in $(\ref{quantityR})$ $R = \infty$. Recall that $\nu_k = \mu_{k+1}$, so $R = \sum_{k=1}^\infty k\nu^{-1}_k = \\ =\sigma^{-2}\sum_{k=1}^\infty k\lambda^{-2\theta}_k$. Since $\lambda > 1$ and $\theta \geq 1$, $R < \infty$, and the result follows.
\cbdu

\begin{Theorem}(Theorem 3.2.3 in \cite{And})\label{forward_equ_uni_thm}
Let Q be the q-matrix as in Theorem \ref{backward_equ_uni_thm}, and assume that $\nu_n > 0$ for all $n \geq 1$. Define

\begin{equation}\label{quantityS}
    S = \sum_{k=1}^\infty \frac{1}{\mu_{k+1}}\Big(1+\frac{\nu_k}{\mu_k}+\frac{\nu_k\nu_{k-1}}{\mu_k\mu_{k-1}}+...+\frac{\nu_k\nu_{k-1}...\nu_2\nu_1}{\mu_k\mu_{k-1}...\mu_2\mu_1} \Big)
\end{equation}
Then the equation $yQ=\lambda y$, $y \in l^+_1$, has unique solution y = 0 if and only if $S = \infty$. Thus, if the minimal solution is dishonest, it is the unique solution of the forward equations if and only if $S = \infty$.
\end{Theorem}

\begin{Corollary}\label{uni-forwardequ}
The forward equation of $(\ref{birth_death_process})$ admits a unique solution.
\end{Corollary}
\pf
First notice that $\sum_{j=1}^\infty e_j(t) \leq C$, where C is an absolute constant, for any $t>0$. Then by Theorem $\ref{forward_equ_uni_thm}$ it suffices to show that $S = \infty$. Indeed, \\
\begin{equation}\notag
    S = \sigma^{-2}\sum_{k=1}^\infty \frac{1}{\lambda^{2k\theta}}\Big (1+\lambda^{2\theta}+\lambda^{4\theta}+...+\lambda^{2k\theta} \Big) \geq \sigma^{-2}\sum_{k=1}^\infty 1 = \infty 
\end{equation}
\cbdu
\begin{remark}
It is worth noticing that although the forward and backward equation have the exactly same form, their solutions are from different spaces. A solution to the forward equation will need to be in $l_1^+$, whereas that to the backward equation is only required to be $l_\infty^+$ with $||y||_{l^+_\infty} \leq 1$. 
\end{remark}
\subsection{The law of the escape time}
In order to prove the anomalous dissipation, we will need to study the law of $\tau$, the time for the minimal process to escape to infinity. We will start from the following lemma:\\
\begin{Lemma}\label{number_visit_k}
 Suppose that $\sum_{j=1}^\infty \lambda^{-2\theta}_j < \infty$ and that the minimal process $\eta_t$ starts from 1. Then for any $k \geq 1$, the number of times $\eta_t$ visits the state $k$ is a geometric random variable with mean equals to $\big ( \lambda^2_{k-1}+\lambda^2_{k} \big ) \sum_{j=k}^\infty \lambda^{-2\theta}_j$
\end{Lemma}
\pf
From the ideas of \cite{Feller}, let $q_{i,j}$ be the transition probability from state $i$ to $j$ of the discrete time Markov chain embedded in continuous-time minimal process and let $\psi^{(k)} = \{\psi_n^{(k)}\}_{n>k}$ denote the probabilities that the process starting from state $n > k$ which will never get back to $k$. Then $\psi^{(k)}$ is the maximal solution of 
\begin{equation}\label{transition_equ}
    x_n = \sum_{j>k}q_{n,j}x_j, \qquad n > k
\end{equation}
such that $0\leq x_n \leq 1$ for all $n$. Note that 
\begin{equation}\notag
    q_{i,j} = \delta_{i,i-1}\Pi_{i,i-1}+\delta_{i,i+1}\Pi_{i,i+1}
\end{equation}\\
where $\delta_{i,j}$ is the kronecker delta. The system of equations $(\ref{transition_equ})$ can be written as
\begin{equation*}
    x_n = \frac{\mu_n}{\mu_n+\nu_n}x_{n-1}+\frac{\nu_n}{\mu_n+\nu_n}x_{n+1}, \qquad  n \geq k+1,
\end{equation*}\\
and it does no harm to assume that $x_k = 0$. Thereafter we get
\begin{equation*}
    x_{n+1}-x_n = \frac{\mu_n}{\nu_n}(x_n-x_{n-1}) = \lambda_1^{-2\theta}(x_n-x_{n-1})
\end{equation*}
and hence inductively, for $n \geq k$
\begin{equation*}
    x_{n+1}-x_n = \lambda^{-2\theta}_{n-k} x_{k+1}
\end{equation*}
which finally yields
\begin{equation*}
    x_n = x_{k+1}\Big ( \lambda^{2\theta}_k \sum^{n-1}_{j=k} \lambda^{-2\theta}_j \Big )
\end{equation*}
The maximum solution is obtained by taking $x_{k+1}$ such that $\lim x_n = 1$, i.e.
\begin{equation*}
    \psi^{(k)}_{k+1} = \Big ( \lambda^{2\theta}_k \sum^{\infty}_{j=k} \lambda^{-2\theta}_j \Big )^{-1},
\end{equation*}
from which we can see that the chain is transient.\\\\
If the minimal process starts from 1, then for any state $k > 1$, it will visit $k$ at least once. The probability that it visits state $k$ for the last time is then given by $q_{k,k+1}\psi^{(k)}_{k+1}$ each time when the process arrives $k$. By the strong Markov property, the total number of visits to $k$ is governed by a geometric random variable $N_k$ with probability of success $q_{k,k+1}\psi^{(k)}_{k+1}$, and whose mean is equal to 
\begin{equation*}
    \mathbb{E}(N_k) = (q_{k,k+1}\psi^{(k)}_{k+1})^{-1} = ( \lambda^2_{k-1}+\lambda^2_{k} \big ) \sum_{j=k}^\infty \lambda^{-2\theta}_j
\end{equation*}
\cbdu
Our main proposition in this section requires us to introduce the following random variable.
\begin{Definition}
For any $n \geq 1$, the total time that the minimal process spends in the state $n$ is denoted by 
\begin{equation*}
    T_n := \mathcal{L} \{t \geq 0: \eta_t = n \}
\end{equation*}
and so $\tau = \sum_{n=0}^\infty T_n$.
\end{Definition}
\begin{Proposition}\label{law_escapeT}
Suppose $r_\infty = \sum^{\infty}_{n=1}n\lambda^{-2\theta}_n < \infty $ and that the minimal process $\eta_t$ start from 1. Then for all $n \geq 1$, $T_n$ is an exponential random variable with mean $\sigma^2r_n, \text{where } r_n:= \sum^\infty_{j=n}\lambda^{-2\theta}_{j}$. In particular, we have 
\begin{equation*}
    \mathbb{E}^{\mathbb{P}}(\tau) = \sum^{\infty}_{n=1}r_n = \sum^{\infty}_{n=1}n\lambda^{-2\theta}_n = r_\infty
\end{equation*}
Moreover, there exists $\alpha > 0$ such that for all $t$, 
\begin{equation*}
    \exp \Big ({-\frac{t}{\sigma^2r_1}}\Big ) \leq \mathbb{P}(\tau > t) \leq \exp \Big(-\frac{t}{\sigma^2r_\infty}+\alpha \Big )
\end{equation*}
\end{Proposition}
\pf
Let $f_{w_n}$ denote the probability distribution function of the waiting time of the process at the state $n$. We know that the waiting time is an exponential distribution with rate $\mu_n+\nu_n$. The probability distribution function for $T_n$ is 
\begin{equation*}
    f_{T_n} = \sum^\infty_{j=1} f^{*(j)}_{w_n} \mathbb {P}\{N_n = j \}
\end{equation*}
where $f^{*(j)}_{w_n}$ is the $j$-fold convolution of $j$ i.i.d exponential distributions, which corresponds to the sum of j copies of i.i.d exponential random variable. Since this sum is a gamma distribution (see for example \cite{PMKS}, 1.4.4), the density function is 
\begin{equation*}
    f^{*(j)}_{w_n}(x) = \frac{\mu_n+\nu_n}{\Gamma(j)}\Big((\mu_n+\nu_n)x\Big)^{j-1}e^{-(\mu_n+\nu_n)x}
\end{equation*}
Therefore by lemma $\ref{number_visit_k}$,\\
\begin{align*}
    f_{T_n} &= \sum^\infty_{j=1}\frac{\mu_n+\nu_n}{\Gamma(j)}\Big((\mu_n+\nu_n)x\Big)^{j-1}e^{-(\mu_n+\nu_n)x} (1-q_{n,n+1}\psi^{(n)}_{n+1})^{j-1}q_{n,n+1}\psi^{(n)}_{n+1}
\end{align*}
 \begin{equation*}
   \Rightarrow \ \ \ \ f_{T_n} = \frac{1}{\sigma^2r_n}e^{-\frac{1}{\sigma^2r_n}x} 
\end{equation*}
This indicates that $T_n$ is an exponential random variable with mean equals to $\sigma^2r_n$.  \\\\
The lower bound is directly from the fact that $\tau = \sum^\infty_{n=1}T_n$,
\begin{equation*}
    \mathbb{P}\{\tau > t \} \geq \mathbb{P} \{T_1 > t \} = e^{-\frac{t}{\sigma^2r_1}}
\end{equation*}
As for the upper bound, we first define $A$ and $\alpha$ respectively as follows:
\begin{align*}
    A &:= -\sum^\infty_{n=1} \sigma^2 r_n \ln (\sigma^2r_n)\\
 \alpha &:= \frac{A}{\sigma^2r_\infty}+\ln (\sigma^2 r_\infty) = -\sum^\infty_{n=1}\frac{r_n}{r_\infty}\ln{\Big(\frac{r_n}{r_\infty}\Big)}
\end{align*}
Then we ought to show that 
\begin{equation*}
    \mathbb{P}\{\tau > t \} \leq \sigma^2 r_\infty \exp{\Big( \frac{A-t}{\sigma^2r_\infty}\Big)}
\end{equation*}
We only need to consider for $t > \sigma^2 r_\infty \ln{(\sigma^2r_\infty)} + A$, otherwise the right hand side of the above inequality is greater than or equal to 1 and we are done. Now for any $n\geq 1 $, let 
\begin{equation*}
    \gm_n := -\frac{\sigma^2r_n}{t}\ln{\big(\sigma^2r_n \big)}-\frac{1}{t}\frac{Ar_n}{r_\infty}+\frac{r_n}{r_\infty}
\end{equation*}
so that 
\begin{equation*}
    \exp{\Big(-\frac{\gm_n t}{\sigma^2r_n} \Big)} = \sigma^2 r_n \exp{\Big(\frac{A-t}{\sigma^2r_\infty} \Big)}
\end{equation*}\\
Notice that since $r_n < r_\infty$ and $t > \sigma^2 r_\infty \ln{(\sigma^2r_\infty)} + A$, $\gm_n > 0$, in addition, 
\begin{equation*}
    \sum_{n=1}^\infty \gm_n  = \sum_{n=1}^\infty \Big [- \frac{\sigma^2r_n}{t}\ln{\big( \sigma^2 r_n\big)} +\frac{1}{t}\frac{Ar_n}{r_\infty}\Big]+1 = 1
\end{equation*}
Therefore we get 
\begin{align*}
    \mathbb{P}\{\tau > t \} &\leq \mathbb{P} \Big ( \bigcup^\infty_{n=1}\{T_n > \gm_n t \} \Big ) \leq \sum^\infty_{n=1} \mathbb{P} \{T_n > \gm_n t \} \\
     &= \sum^\infty_{n=1}  \exp{\Big(-\frac{\gm_n t}{\sigma^2r_n} \Big)}= \sum^\infty_{n=1}\sigma^2 r_n \exp{\Big(\frac{A-t}{\sigma^2r_\infty} \Big)}\\
     &=\sigma^2 r_\infty \exp{\Big(\frac{A-t}{\sigma^2r_\infty} \Big)} = \exp{(-\frac{t}{\sigma^2 r_\infty}+\alpha)}
\end{align*}
\cbdu
The next lemma shows that either $\tau = \infty$ a.s. or $\mathbb{P} \{ \tau \in [T_1,T_2]\} > 0 $ for any $T_1 < T_2$.   
\begin{Lemma}\label{lemma1_escapeT}
Suppose $\mathbb{P}\{\tau > T \} < 1$ for some $T$. Then for any $T_1 < T_2$ we have that
$\mathbb{P}\{\tau > T_1\} > \mathbb{P}\{\tau > T_2\}$. In particular, $\mathbb{P}\{\tau > t \} < 1$ for all $t>0$.
\end{Lemma}
\pf
First of all, from the definition of $\tau$: 
\begin{equation*}
    \mathbb{P} \{\tau > t\} = \sum^\infty_{i=1} \mathbb{P}\{\eta_t = i\}
\end{equation*}
In addition, we denote
\begin{equation*}
    \mathbb{P}\{\tau > t | \eta_0 = m \} = \sum^\infty_{n=1}\mathbb{P}\{\eta_t = n|\eta_0 = m \}:=\sum_{n=1}^\infty p_{m,n}(t)
\end{equation*}
and also recall that  
\begin{equation*}
    \mathbb{P} \{ \eta_t = j\} = e_j(t) 
\end{equation*}\\
By Chapman-Kolmogorov equation (see \cite{Feller}), that is $p_{m,n}(s+t) = \sum_{j=1}^\infty p_{m,j}(s)p_{j,n}(t)$, we obtain the following: 
\begin{align*}
    \mathbb{P}\{\tau > t \} &= \sum^\infty_{m=1}\mathbb{P}\{\tau > t | \eta_0 = m \}\mathbb{P}\{\eta_0 = m\}\\
    &= \sum_{m=1}^\infty \sum_{n=1}^\infty p_{m,n}(t) e_m(0)\\
    &= \sum_{m=1}^\infty \sum_{n=1}^\infty \sum_{j=1}^\infty p_{m,j}(t-s)p_{j,n}(s) e_m(0)\\
    &= \sum_{n=1}^\infty \sum_{j=1}^\infty e_j(t-s)p_{j,n}(s) = \sum_{j=1}^\infty \mathbb{P}\{\tau > s | \eta_0 = j \} \mathbb{P} \{\eta_{t-s} = j\}\\
    & \leq \sum_{j=1}^\infty \mathbb{P} \{\eta_{t-s} = j \} = \mathbb{P} \{\tau > t-s \}
\end{align*}
Now suppose that $\mathbb{P}\{\tau > T_1\} = \mathbb{P}\{\tau > T_2\}$ for some $T_1 < T_2$ then if we take $s = T_2 - T_1$ and $t = T_2$, the above inequality becomes equality. This implies that $\mathbb{P}\{\tau > T_2 -T_1 | \eta_0 = j \} = 1$ for all $j$., hence $\mathbb{P}\{\tau > T_2 -T_1\} = 1$. On the other hand we know that $\mathbb{P} \{\tau > 0 \} = 1$. Since  $\mathbb{P}\{\tau > t \}$ is non-increasing, we conclude that $\mathbb{P}\{\tau > t \} = 1$ on $[0,T_2 - T_1]$. Repeat this procedure for $t = 2(T_2-T_1)$ and $s = T_2-T_1$, and inductively we conclude that $\mathbb{P} \{\tau > 0 \} \equiv 1$, contradicting to the assumption. 
\cbdu
We conclude this section with the following lemma, which states that the probability of no explosion is greater if the process $\eta_t$ starts from 1. 
\begin{Lemma}\label{lemma2_escapeT}
$\mathbb{P} \{\tau > t \} \leq \mathbb{P} \{\tau > t \ |\ \eta_0 = 1 \}$
\end{Lemma}
\pf
Firstly, let us define the random variable $H_k$ as the first time when $\eta_t$ hitting the state $k$, where $k > 1$. Since the process is transient, the distribution of $H_k$ satisfies that $\mathbb{P}(H_k \leq s) := F(s) \to 1$ as $s \to \infty$. By the strong Markov Property and lemma $\ref{lemma1_escapeT}$, 
\begin{align*}
    \mathbb{P} \{\tau > t \ |\  \eta_0 =1\} &= \int_{0}^\infty \mathbb{P} \{\tau > t \ |\  \eta_0 =1,\ H_k=s \} dF(s) \\
    &= \int_{0}^\infty \mathbb{P} \{\tau > t \ |\  \eta_0 =1,\ \eta_s = k \} dF(s)\\
    &= \int_{0}^\infty \mathbb{P} \{\tau > t \ |\ \eta_s = k \}dF(s) \\
    &= \int_{0}^\infty \mathbb{P} \{\tau > t-s \ |\ \eta_0 = k \}dF(s)\\
    &\geq \int_{0}^\infty \mathbb{P} \{\tau > t \ |\ \eta_0 = k \}dF(s) = \  \mathbb{P} \{\tau > t \ |\ \eta_0 = k \}
\end{align*}\\
The above argument applies for all $k > 1$ , hence we deduce that \\
\begin{align*}
    \mathbb{P} \{\tau > t \} &= \sum_{k = 1}^\infty \mathbb{P} \{\tau > t \ |\  \eta_0 =k\} \mathbb{P} \{\eta_0 = k \}\\
    &\leq \sum_{k = 1}^\infty \mathbb{P} \{\tau > t \ |\  \eta_0 =1\} \mathbb{P} \{\eta_0 = k \} \\
    & =  \mathbb{P} \{\tau > t \ |\  \eta_0 =1\} \sum_{k = 1}^\infty e_k(0)\\
    & =  \mathbb{P} \{\tau > t \ |\  \eta_0 =1\}
\end{align*}
\cbdu
\subsection{Decay of the energy $e(t)$} In this section, we will prove that the energy $e(t)$ goes to 0 as $t \to \infty$ both in $L^1$ and almost surely under $\tilde{Q}$. Moreover, we will show that given the initial energy profile $e(0)$ sufficiently small, $e(t)$ decay similarly under the measure $Q$. 
\begin{Proposition}\label{decay in L^1 under Q_tilde}
Suppose $r_\infty := \sum_{n=1}^\infty n\lambda^{-2}_n < \infty$. Let $X(t) = (P(t),M(t))$ be the solution in the Theorem $\ref{thm-main}$ and $\Upsilon(t) := \frac{1}{2}\sum_{n=1}^\infty \big( P^2_n(t)+M^2_n(t) \big)$ the energy at time $t$. Then 
\begin{equation}\label{L1_conv}
    \lim_{t \to \infty } E^{\tilde {Q}} [ \Upsilon(t) ] = 0 
\end{equation}
\end{Proposition}
\pf
From the definition, 
\begin{equation*}
    E^{\tilde {Q}} [ \Upsilon(t) ] = \frac{1}{2} \sum_{n=1}^\infty E^{\tilde {Q}} [ P^2_n(t)+M^2_n(t) ] = \Upsilon(0) \sum_{n=1}^\infty e_n(t)
\end{equation*}
so by Corollary $\ref{uni-forwardequ}$ and Proposition $\ref{law_escapeT}$, we have for any $t > 0$
\begin{equation*}
    \Upsilon(0) \sum_{n=1}^\infty e_n(t) = \Upsilon(0) \sum_{n=1}^\infty \mathbb{P} \{\eta_t = n\} = \Upsilon(0) \mathbb{P} \{\tau > t \} \leq \Upsilon(0)\exp \Big(-\frac{t}{\sigma^2r_\infty}+\alpha \Big )
\end{equation*}\\
whence,
\begin{equation}\label{L_1 bound of energy under Q_tilde}
     E^{\tilde {Q}} [ \Upsilon(t) ] \leq \Upsilon(0)\exp \Big(-\frac{t}{\sigma^2r_\infty}+\alpha \Big )
\end{equation}
for some $\alpha > 0$. Sending t to infinity we obtain $(\ref{L1_conv})$. 
\cbdu
\vspace{0.3cm}
In order to show that $\Upsilon$ goes to 0 a.s. under $\tilde{Q}$, we need the following lemma:
\begin{Lemma}\label{Energy_monotone_a.s}
Let $X(t)$ and $\Upsilon(t)$ as above. Then for all $t \geq s \geq 0 $ we have that \\
\begin{equation}\notag
    \tilde{Q} \Big ( \Upsilon (t) \leq \Upsilon(s) \Big) = 1   
\end{equation}
\end{Lemma}
\pf
    For any $s \geq 0$, from the Theorem \ref{thm-G}, with the initial condition $X(s) = (P_j(s), M_j(s))$ there exists a unique solution $\Hat{X}$ such that 
    \begin{equation}
         \Tilde{Q} \Big ( \sum_{j=1}^\infty \Hat{X}^2_j (t)  \leq \sum_{j=1}^\infty X^2_j(s)  \Big) = 1
    \end{equation}
    for every $t \geq s$. On the other hand, from the uniqueness $\Hat{X}(t) = X(t)$, hence the lemma follows.
\cbdu

\begin{Proposition}\label{decay_a.s._under Q_tilde} 
With the same assumption as Proposition $\ref{decay in L^1 under Q_tilde}$, the total energy $\Upsilon(t)$ has at least an exponential decay under $\tilde{Q}$: \\
\begin{equation}
    \limsup_{t \to \infty}{\frac{1}{t}\ln \Upsilon(t)} \leq -\frac{1}{\sigma^2 r_\infty} \quad \tilde{Q} - a.s.
\end{equation}
\end{Proposition}
\pf
For any $\epsilon > 0$ , we define 
\begin{equation}\notag
    \beta := -\frac{1}{\sigma^2 r_\infty} + \epsilon
\end{equation}
then for any positive integer $n$, by proposition $\ref{decay in L^1 under Q_tilde}$, \\
\begin{align*}
    \tilde{Q} \Big (n^{-1}\ln{\Upsilon (n)} > \b \Big ) \leq e^{-\b n}\mathbb{E}^{\tilde{Q}}\Big (\Upsilon(n) \Big) \leq C_1 e^{-\e n}
\end{align*}
where $C_1  = \Upsilon(0) e^\a$, which is independent of $n$. Therefore, by Borel-Cantelli lemma, there exists a null set $W$ which satisfies the following: $\forall \omega \in W^c$, $\exists N_0(\omega)$ such that if $n\geq N_0(\om)$ then $\Upsilon(n,\om) \leq e^{-\b n}$. Let 
\begin{equation}\notag
    C_2(\om) := \sup_{n<N_0}\Upsilon(n,\omega)e^{\b n}
\end{equation}
From lemma $\ref{Energy_monotone_a.s}$ , there is another null set $\tilde{W}$ such that 
\begin{equation}\notag
    \Upsilon(t,\om) \leq \Upsilon (\lfloor t \rfloor,\om), \quad \forall \om \in \Tilde{W}^c, \ \forall t \in [0,\infty)
\end{equation}
Hence for all $\om \in W^c\cap \Tilde{W}^c$ , we have that 
\begin{equation}\notag
    \Upsilon(t,\om) \leq \Upsilon(\lfloor t \rfloor,\om) \leq C_2(\om)e^{-\b \lfloor t \rfloor} \leq 2 C_2(\om)e^{-\b t} 
\end{equation}
and then sending $\epsilon$ to 0 completes the proof.
\cbdu

It should be well noticed that the similar decaying properties need not to be held under the measure $Q$, because the constant $C_2(\om)$ is not $F_t$ -measurable but $F_\infty$-measurable. The following proposition provides a sufficient condition for the equivalence of $Q$ and $\Tilde{Q}$, hence the decaying of $\Upsilon (t)$ $a.s.$ and in $L^1$ follows.

\begin{Proposition}\label{Equivalent_measure}
Let $X$ be the weak solution obtained from the Theorem \ref{thm-main}. Assume that $r_\infty = \sum_{n=1}^\infty n\lambda^{-2\theta}_n < \infty$ and that $\Upsilon(0) < \frac{1}{\sigma^2r_\infty}$ then 
\begin{equation}\notag
    \mathbb{E}^{\Tilde{Q}} \big(e^{\int_0^\infty \Upsilon(t)dt}   \big) <\infty 
\end{equation}
hence $Q$ and $\Tilde{Q}$ are equivalent on $F_\infty$.
\end{Proposition}
\pf
From Proposition \ref{decay in L^1 under Q_tilde} and Lemma $\ref{Energy_monotone_a.s}$, we have 
\begin{align*}
 E^{\tilde {Q}} [ \Upsilon(t) ] &\leq \Upsilon(0)\exp \Big(-\frac{t}{\sigma^2r_\infty}+\alpha \Big ), \quad \forall t \geq 0\\
 \Upsilon(t) &\leq \Upsilon(0), \quad \Tilde{Q}-a.s. 
\end{align*}
Define the non-negative random variable $Z := \int_0^\infty \Upsilon(t)dt$. For any $z \geq 0$, we deduce the following estimation
\begin{align*}
    z\Tilde{Q}\Big (Z>z \Big) &= z\mathbb{E}^{\tilde{Q}} \big [ \mathds{1}_{Z>z} \big] \leq \mathbb{E}^{\Tilde{Q}} \big( Z \cdot \mathds{1}_{Z>z}\big) \\
    &= \int_0^\infty \mathbb{E}^{\Tilde{Q}} \big( \Upsilon(t) \cdot \mathds{1}_{Z>z}\big) dt\\
    &\leq \int_0^\infty \min \Big( \mathbb{E}^{\Tilde{Q}} \big( \Upsilon(t) \big), \Upsilon(0)\Tilde{Q}(Z>z) \Big) dt\\
    &\leq \Upsilon(0)\int_0^\infty \min \Big( \exp{\big (-\frac{t}{\sigma^2 r_\infty}+\alpha \big)}, \Tilde{Q}(Z>z) \Big) dt\\
    &=\Upsilon(0) \int_0^s \Tilde{Q}(Z>z)dt   + \Upsilon(0)\int_s^\infty \exp{\big( -\frac{t}{\sigma^2 r_\infty}+\alpha \big)}dt 
\end{align*}
where $s$ is such that $\Tilde{Q}(Z>z) = \exp{\big( -\frac{s}{\sigma^2 r_\infty}+\alpha \big)}$, and from which yields 
\begin{align*}
    z\Tilde{Q}(Z > z) &\leq \Upsilon(0) \Big ( s\Tilde{Q}(Z > z)+ \sigma^2\ r_\infty \exp{\big( -\frac{s}{\sigma^2 r_\infty}+\alpha \big)}\Big)\\
   \Rightarrow z &\leq \Upsilon(0)\big(s + \sigma^2 r_\infty \big)\\
   \Rightarrow s &\geq \frac{z}{\Upsilon(0)} - \sigma^2 r_\infty
\end{align*}
We can assume that $\tilde{Q}(Z > z) > 0$ for any $z > 0$, since otherwise $Z$ is bounded a.s. and the result is automatic. Hereafter, 
\begin{align*}
    \Tilde{Q} \big ( Z > z \big ) = \exp{\big( -\frac{s}{\sigma^2 r_\infty}+\alpha \big)} &\leq \exp{\big( -\frac{z}{\sigma^2 r_\infty \Upsilon(0)}+\alpha+1 \big)}\\
    \Rightarrow \Tilde{Q} \big ( e^Z > w \big ) &\leq w^{-\frac{1}{\sigma^2 r_\infty \Upsilon(0)}}e^{\alpha + 1}
\end{align*}
These lead to the result : 
\begin{equation*}
    \mathbb{E}^{\tilde{Q}} (e^Z)  = \int_0^\infty \Tilde{Q}(e^Z > w)dw \leq 1 + e^{\a +1}\int_1^\infty w^{-\frac{1}{\sigma^2 r_\infty \Upsilon(0)}}dw < \infty
\end{equation*}
since $\sigma^2 r_\infty \Upsilon(0) < 1$.
\cbdu
\vspace{0.5cm}
\begin{Corollary}\label{exp_decay a.s. under Q}
With the same assumption as Proposition \ref{Equivalent_measure}, the energy function $\Upsilon(t)$ decay at least exponentially fast under the measure $Q$: 
\begin{equation*}
    \limsup_{t \to \infty}{\frac{1}{t}\ln \Upsilon(t)} \leq -\frac{1}{\sigma^2 r_\infty} \quad Q- a.s.
\end{equation*}
\end{Corollary}

\pf
Directly from Proposition \ref{Equivalent_measure}, the measure Q and $\Tilde{Q}$ are equivalent. So by Proposition \ref{decay_a.s._under Q_tilde} the similar result holds Q - a.s.
\cbdu
The next proposition shows further that if the energy is bounded by $\frac{1}{ r_\infty}$, then a similar decay of energy profile in $L^1$ under the measure $Q$ is expected.
\vspace{0.5cm}
\begin{Proposition}\label{decay in L^1 under Q}
Let $C$ be the constant in Definition \ref{def-bounded}. If $C < \frac{1}{r_\infty}$, then 
\begin{equation*}
    \limsup_{t \to \infty}{\frac{1}{t}\ln{\mathbb{E}^Q\big(\Upsilon(t) \big)}} \leq -\frac{1}{\sigma^2 r_\infty} \big (1-\sqrt{Cr_\infty} \big)^2
\end{equation*}
\end{Proposition}
\pf
Denote $D_t$ by 
\begin{equation*}
    D_t := \exp{\Big(-\big(L_{1,t}+L_{2,t} \big)+ \frac{1}{2}\big([L_1,L_1]_t+[L_2,L_2]_t \big)\Big)}
\end{equation*}
where $L_{1,t}$ and $L_{2,t}$ were given in Lemma \ref{le-ex-Ito}. Since the solution is in $L^\infty(\Omega\times[0,T]; l^2\times l^2)$ under $Q$, we have that 
\begin{equation*}
    \exp{\Big( \frac{\gm}{2}\big( [L_1,L_1]_t+[L_2,L_2]_t \big ) \Big ) } \leq \exp{(\frac{C\gm t}{\sigma^2})} \quad Q-a.s.
\end{equation*}
for any $\gm > 0$. Here C is the constant such that $\Upsilon(t) \leq C$ Q-a.s. for all $t \geq 0$. Use Lemma \ref{Energy_monotone_a.s} and (\ref{L_1 bound of energy under Q_tilde}) we obtain
\begin{align*}
    \mathbb{E}^Q \big ( \Upsilon(t) \big) &= \mathbb{E}^{\Tilde{Q}} \big ( \Upsilon(t) D_t\big) \leq  \bigg ( \mathbb{E}^{\Tilde{Q}} \big ( D^p_t\big) \bigg)^\frac{1}{p} \bigg ( \mathbb{E}^{\Tilde{Q}} \big ( \Upsilon^q(t) \big) \bigg)^\frac{1}{q}\\
    &\leq  \bigg ( \mathbb{E}^{\Tilde{Q}} \big ( D^p_t\big) \bigg)^\frac{1}{p} \bigg ( \mathbb{E}^{\Tilde{Q}} \big ( \Upsilon(t) \Upsilon(0)^{q-1} \big) \bigg)^\frac{1}{q} \\
    &\leq  \Upsilon(0)^{1-\frac{1}{q}}\bigg ( \mathbb{E}^{\Tilde{Q}} \big ( D^p_t\big) \bigg)^\frac{1}{p} \bigg ( \Upsilon(0)\exp \Big(-\frac{t}{\sigma^2r_\infty}+\alpha \Big )\bigg)^\frac{1}{q} \\
    &\leq \Upsilon(0) \bigg ( \mathbb{E}^{\Tilde{Q}} \big ( D^p_t\big) \bigg)^\frac{1}{p} \exp \Big(-\frac{t}{q\sigma^2r_\infty}+\frac{\a}{q} \Big )
\end{align*}
where $p, q > 1 $ such that $\frac{1}{p}+ \frac{1}{q} = 1 $. The p-moment of $E_t$ under $\Tilde{Q}$ can be estimated by 
\begin{align*}
     \mathbb{E}^{\Tilde{Q}} \big ( D^p_t\big) &=  \mathbb{E}^{Q} \big ( D^{p-1}_t\big) \\
     &= \mathbb{E}^{Q} \Big [ \exp{\big(-(p-1)(L_{1,t}+L_{2,t} ) +\frac{p-1}{2}([L_1,L_1]_t+[L_2,L_2]_t ) \big) } \Big]\\
     &= \mathbb{E}^{Q} \Big[ E_t \exp{\Big(\frac{p(p-1)}{2} ([L_1,L_1]_t+[L_2,L_2]_t \Big)} \Big]
\end{align*}
where 
\begin{equation*}
    E_t := \exp{\big(-(p-1)(L_{1,t}+L_{2,t} ) -\frac{(p-1)^2}{2}([L_1,L_1]_t+[L_2,L_2]_t ) \big) }
\end{equation*}
Since the solution is exponentially integrable, by Novikov condition, $E_t$ is a martingale: 
\begin{equation*}
    \mathbb{E}^Q (E_t) = \mathbb{E}^Q (E_0)= 1
\end{equation*}
and so 
\begin{align*}
    \mathbb{E}^{\Tilde{Q}} \big ( D^p_t\big) &=\mathbb{E}^{Q} \Big[ E_t \exp{\Big(\frac{p(p-1)}{2} ([L_1,L_1]_t+[L_2,L_2]_t \Big)} \Big] \\
     &\leq e^{C\frac{p(p-1)}{\sigma^2})t} \mathbb{E}^{Q} \Big[ E_t \Big] = e^{C\frac{p(p-1)}{\sigma^2}t} 
\end{align*}
In conclusion,
\begin{equation*}
    \mathbb{E}^Q \big ( \Upsilon(t) \big) \leq \Upsilon(0)\exp \Big(\frac{C(p-1)t}{\sigma^2}-\frac{t}{q\sigma^2r_\infty}+\frac{\a}{q} \Big )
\end{equation*}
Next, we seek for the p minimizing the right hand side. It turns out for a fixed $t > \frac{\a \sigma^2}{1/r_\infty - C}$, where we need  this $p = \sqrt{\frac{1/r_\infty - \a \sigma^2/t}{C}} > 1$. Computation shows that with such $p$,
\begin{align*}
    \mathbb{E}^Q \big ( \Upsilon(t) \big) &\leq \Upsilon(0)\exp \Big(\frac{C(p-1)t}{\sigma^2}-\frac{Cp(p-1)t}{\sigma^2} \Big )\\
    &=\Upsilon(0)\exp \Big(-\frac{C(p-1)^2t}{\sigma^2} \Big )
\end{align*}
which implies that 
\begin{align*}
    \frac{1}{t}\ln{\mathbb{E}^Q \big ( \Upsilon(t) \big)} &\leq \frac{1}{t}\Upsilon(0) - \frac{C(p-1)^2}{\sigma^2}\\
    &= \frac{1}{t}\Upsilon(0) - \frac{1}{\sigma^2 r_\infty} \Big(\sqrt{1-\frac{\a r_\infty C \sigma^2}{t}} -\sqrt{C r_\infty}\Big )^2 \\
\Rightarrow \limsup_{t \to \infty}\frac{1}{t}\ln{\mathbb{E}^Q \big ( \Upsilon(t) \big)} &\leq -\frac{1}{\sigma^2 r_\infty} \big (1-\sqrt{Cr_\infty} \big)^2
\end{align*}
\cbdu
We are ready to prove the Theorem $\ref{anomalous_dissipation}$ now.
\subsection{Proof of Theorem \ref{anomalous_dissipation}} Since $r_\infty = \sum_{n=1}^\infty n\lambda^{-2\theta}_n < \infty$, by Proposition \ref{law_escapeT} and Lemma \ref{lemma2_escapeT}, there exist t > 0 such that $\mathbb{P} \{ \tau > t \} < 1$. Lemma \ref{lemma1_escapeT} further implies that $\mathbb{P} \{ \tau > t \} < 1$ for all $t > 0$. From the proof of Proposition \ref{decay in L^1 under Q_tilde}, 
\begin{equation*}
    \mathbb{E}^{\tilde{Q}} \big( \Upsilon(t)  = \Upsilon(0)\big)\mathbb{P} \{\tau > t \} < \Upsilon(0)
\end{equation*}
and together with Lemma \ref{Energy_monotone_a.s} we obtain that 
\begin{equation*}
    \Tilde{Q} \big (\Upsilon(t) = \Upsilon(0) \big) <1 , \quad \forall t>0
\end{equation*}
In addition, Proposition \ref{decay in L^1 under Q_tilde} also yields that 
\begin{equation*}
    \Tilde{Q} \big(\Upsilon(t) > \epsilon \big) \leq \frac{\mathbb{E}^{\tilde{Q}} \big( \Upsilon(t) \big)}{\epsilon} \leq \frac{1}{\epsilon}\Upsilon(0) \exp{\big(-\frac{t}{\sigma^2r_\infty}+\a \big)} < 1
\end{equation*}
for $t$ sufficiently large. The same results hold under the measure $Q$ as two measures are equivalent on $F_t$ for all $t > 0$. The last statement is proven by Corollary \ref{exp_decay a.s. under Q} and Proposition \ref{decay in L^1 under Q}.

\subsection{Lack of solution with regularity}
Theorem \ref{anomalous_dissipation} indicates the lack of regularity of the solution obtained from Theorem \ref{thm-main}. This can be shown by the following:
\begin{Definition}
We define the space $H$ by
\begin{align*}
    H := \bigg \{ Z = (u,v) \in l^2\times l^2: \sum_{j=1}^\infty \lambda^{2\theta}_j (u^2_n+ v^2_n) < \infty   \bigg \}
\end{align*}
H is a Hilbert space with norm $||Z||^2_H = \sum_{j=1}^\infty \lambda^{2\theta}_j (u^2_n+ v^2_n)$
\end{Definition}
\begin{Proposition}
The solution $X(t)= (P(t),M(t))$ obtained by Theorem \ref{thm-main} satisfies 
\begin{equation*}
    Q \Big ( \int_0^T ||X(t)||^2_H dt = \infty \Big) > 0
\end{equation*}
\end{Proposition}

\pf
We will prove the statement by contradiction. Suppose that 
\begin{equation*}
   Q \Big ( \int_0^T ||X(t)||^2_H dt < \infty \Big) = 1
\end{equation*}
The Ito's form (\ref{sys-Ito}) of the system yields 
\begin{align*}
    \frac{1}{2} d \bigg( \sum_{j=1}^N \Big (P^2_j(t) + M^2_j(t) \Big) \bigg) = -\lambda^\theta_N M_N P_N P_{N+1} dt  -\lambda^\theta_N P_N M_N M_{N+1} dt\\
    -\sigma \lambda^\theta_N P_N P_{N+1} dW^p_N  -\sigma M_N M_{N+1} dW^m_N -\frac{\sigma^2}{2}\sum_{j=1}^N \big(\lambda^{2\theta}_{j-1} + \lambda^{2\theta}_{j}\big) (P^2_j+M^2_j)dt\\
    +\frac{\sigma^2}{2}\sum_{j=1}^N \Big(\lambda^{2\theta}_{j-1}P^2_{j-1} + \lambda^{2\theta}_{j}P^2_{j+1}+\lambda^{2\theta}_{j-1}M^2_{j-1} + \lambda^{2\theta}_{j}M^2_{j+1} \Big)dt \hspace{2.6cm}
\end{align*}\\
and if we denote $X(0)= (p,m) \in l^2 \times l^2$ then $\forall t \leq T$
\begin{align*}
    \sum_{j=1}^N \Big (P^2_j(t) + M^2_j(t) \Big) - \sum_{j=1}^N \big (p^2_j + m^2_j \big) = -\int_0^t 2\lambda^\theta_N M_N P_N \Big (M_{N+1} + P_{N+1} \Big) ds \\
    -\int_0^t 2\lambda^\theta_N P_N P_{N+1} dW^p_N(s) -\int_0^t 2\lambda^\theta_N M_N M_{N+1} dW^m_N(s) -\sigma^2 \int_0^t \lambda^{2\theta}_N \Big (P^2_N + M^2_N\Big) ds \\
    +\sigma^2 \int_0^t \lambda^{2\theta}_{N} \Big (P^2_{N+1} + M^2_{N+1}\Big) ds \hspace{8cm}
\end{align*}\\
Denote that 
\begin{align*}
    I_N := &-\int_0^t 2\lambda^\theta_N M_N P_N \Big (M_{N+1} + P_{N+1} \Big) ds -\sigma^2 \int_0^t \lambda^{2\theta}_N \Big (P^2_N + M^2_N\Big) ds \\&+\sigma^2 \int_0^t \lambda^{2\theta}_{N} \Big (P^2_{N+1} + M^2_{N+1}\Big) ds
\end{align*}
and that 
\begin{align*}
    J_N := -\int_0^t 2\lambda^\theta_N P_N P_{N+1} dW^p_N(s) -\int_0^t 2\lambda^\theta_N M_N M_{N+1} dW^m_N(s)
\end{align*}
Use the fact that $X(t)\in L^\infty (\Om \times [0,T];l^2 \times l^2)$, there exists a constant $C > 0$ such that 
\begin{equation*}
    |P_{N+1}|+|M_{N+1}| \leq C
\end{equation*}
and also note that $\lambda^{\theta}_N \leq \lambda^{\theta}_{N+1}$. Therefore, we can bound $I_N$ by
\begin{align*}
    |I_N| &\leq C\int_0^t \bigg ( 2\lambda^\theta_N M_N P_N + \sigma^2 \lambda^{2\theta}_N \Big( P^2_N + M^2_N\Big) + \sigma^2 \lambda^{2\theta}_{N+1} \Big( P^2_{N+1} + M^2_{N+1}\Big) \bigg) ds \\
     &\leq C\int_0^t \bigg ( \lambda^{2\theta}_N \Big( M^2_N+P^2_N\Big ) + \sigma^2 \lambda^{2\theta}_N \Big( P^2_N + M^2_N\Big) + \sigma^2 \lambda^{2\theta}_{N+1} \Big( P^2_{N+1} + M^2_{N+1}\Big) \bigg) ds 
\end{align*}\\
By the assumption, $I_N$ goes to 0 as $N \to 0$ a.s. and hence in probability. As for $J_N$, we will illustrate the estimation of the first term. The second one follows by the exact same manner. For any $\epsilon, \delta > 0$, we have that 
\begin{align*}
    Q\bigg( \Big |\int_0^t \lambda^\theta_N P_N P_{N+1} dW^p_N(s) \Big | > \epsilon \bigg) &= Q\bigg( \Big |\int_0^t \lambda^\theta_N P_N P_{N+1} dW^p_N(s) \Big | > \epsilon , Y_t > \delta   \bigg)\\
    &+  Q\bigg( \Big |\int_0^t \lambda^\theta_N P_N P_{N+1} dW^p_N(s) \Big | > \epsilon , Y_t < \delta   \bigg)
\end{align*}
where $Y_t = \int_0^t \lambda^{2\theta}_N P^2_N P^2_{N+1} ds$. In addition, we define a stopping time 
\begin{equation*}
    \tau_\delta := \inf \Big \{s :  \int_0^s  \lambda^{2\theta}_N P^2_N(v) P^2_{N+1}(v) dv = \delta \Big \}
\end{equation*}
Continuing the estimation one derives the following, \\
\begin{align*}
     &Q\bigg( \Big |\int_0^t \lambda^\theta_N P_N P_{N+1} dW^p_N(s) \Big | > \epsilon \bigg) \leq Q \Big( Y_t > \delta \Big ) + \frac{1}{\epsilon^2}\mathbb{E}\bigg( \Big |\int_0^t \lambda^\theta_N P_N P_{N+1} dW^p_N(s) \Big |^2 \mathds{1}_{\{Y_t < \delta\}} \bigg)\\\\
     &\leq Q \Big ( \int_0^t \lambda^{2\theta}_N P^2_N P^2_{N+1} ds > \delta\Big) + \frac{1}{\epsilon^2}\mathbb{E}\bigg( \Big |\int_0^t \lambda^\theta_N P_N P_{N+1} dW^p_N(s) \Big |^2 \mathds{1}_{\{\tau_{\delta} \geq t\}} \bigg)\\\\
   & \leq  Q \Big ( \int_0^t \lambda^{2\theta}_N P^2_N P^2_{N+1} ds > \delta\Big) + \frac{1}{\epsilon^2}\mathbb{E}\bigg( \int_0^{\tau_\delta} \lambda^{2\theta}_N P^2_N P^2_{N+1} ds   \bigg)\\\\
    &\leq  Q \Big ( \int_0^t \lambda^{2\theta}_N P^2_N P^2_{N+1} ds > \delta\Big) + \frac{\delta}{\epsilon^2}
\end{align*}\\\\
Similarly,\\
\begin{equation*}
    Q\bigg( \Big |\int_0^t \lambda^\theta_N M_N M_{N+1} dW^m_N(s) \Big | > \epsilon \bigg) \leq Q \Big ( \int_0^t \lambda^{2\theta}_N M^2_N M^2_{N+1} ds > \delta\Big) + \frac{\delta}{\epsilon^2}
\end{equation*}\\\\
Therefore, $J_N$ goes to 0 in probability as $N \to \infty$. Hence, in probability we have the convergence\\ 
\begin{equation*}
   \frac{1}{2}\sum_{j=1}^N \Big (P^2_j(t) + M^2_j(t) \Big) \xrightarrow[N \to \infty]{}  \Upsilon(0)
\end{equation*}\\
On the other hand we have the almost sure convergence:\\ 
\begin{equation*}
   \frac{1}{2}\sum_{j=1}^N \Big (P^2_j(t) + M^2_j(t) \Big) \xrightarrow[N \to \infty]{}  \Upsilon(t)
\end{equation*}\\
and so \\
\begin{equation*}
    Q \Big(\Upsilon(t)= \Upsilon(0) \Big ) = 1, \qquad \forall t \leq T
\end{equation*}\\
However, from Theorem $\ref{anomalous_dissipation}$ it is a contradiction. We then conclude that \\
\begin{equation*}
    Q \Big ( \int_0^T ||X(t)||^2_H dt = \infty \Big) > 0
\end{equation*}
\cbdu

\section*{Acknowledgement}
M. Dai and Q. Peng are partially supported by NSF Grant DMS--2009422. M. Dai is also grateful for the support of the AMS Centennial Fellowship, and the hospitality of the Institute for Advanced Study and Princeton University where part of the work was completed. 


\bigskip


\end{document}